\documentclass[a4paper,oneside,11pt]{article}

\newtheorem{dfn}{Definicja}[section]
\newtheorem{tw}[dfn]{Theorem}
\newtheorem{prop}[dfn]{Proposition}

\newtheorem{lem}[dfn]{Lemma}
\newtheorem{rem}[dfn]{Remark}
\newtheorem{ex}[dfn]{Example}

\newtheorem{cor}[dfn]{Corollary}

\usepackage{anysize}
\usepackage{verbatim}
\usepackage{array}
\usepackage{enumerate}
\usepackage{amssymb} 
\usepackage{amsmath}
\usepackage{fancyhdr}
\usepackage{euscript}
\usepackage{latexsym}

\author{Micha\l \ Barski \\
\small Faculty of Mathematics and Computer Science, University of
Leipzig, Germany\\
\small Faculty of Mathematics Cardinal Stefan Wyszy\'nski University in Warsaw\\
\small{\small{\it Michal.Barski@math.uni-leipzig.de}}}

\title{\bf Integral representations of risk functions for basket derivatives}

\begin{document}
\baselineskip=1.1\baselineskip \maketitle
\date

\begin{abstract}
The risk minimizing problem
$\mathbf{E}[l((H-X_T^{x,\pi})^{+})]\overset{\pi}{\longrightarrow}\min$
in the multidimensional Black-Scholes framework is studied. Specific
formulas for the minimal risk function and the cost reduction
function for basket derivatives are shown.
Explicit integral representations for the risk functions for $l(x)=x$ and $l(x)=x^p$,
with $p>1$ for digital, quantos, outperformance and spread options
are derived.
\end{abstract}

\noindent
\begin{quote}


\end{quote}

\section{Introduction}
The paper is devoted to the stochastic control problem arising in
the risk analysis of financial markets. Let $H$ be a random variable
representing future random payoff which is traded on the market.
Denote its price  determined by the no arbitrage method by $p(H)$.
If the initial capital $x$ of the writer exceeds $p(H)$ then he is
able to hedge $H$ perfectly, i.e. he can follow some trading
strategy $\pi$ such that the wealth process at the final time is
greater than $H$, i.e.
\begin{gather*}
 P(X_T^{x,\pi}\geq H)=1.
\end{gather*}
If $x<p(H)$ then the above equality fails for each $\pi$ and as a
consequence shortfall risk appears. The aim of the trader is to find
a strategy which is optimal in a sense. Let
$l:[0,+\infty)\longrightarrow[0,+\infty)$ be a loss function which
describes the attitude of the trader to hedging losses. The goal
is to minimize
\begin{gather}\label{shortfall risk}
\mathbf{E}[l((H-X_T^{x,\pi})^{+})].
\end{gather}
This problem was studied with various model settings in many papers.
The ones mentioned below do not form a complete list. Existence of
the optimal strategy for the case when $l(x)=x$ in the context of
complete market with the stock prices modeled by the diffusion
processes was shown in {\cite{CvitanicKaratzas}}. These results were
generalized to incomplete markets in \cite{Cvitanic} where
existence of solution with the use of dual methods was shown.
Existence of the optimal strategy in a general semimartingale model
was shown in \cite{Pham}. In \cite{FL2} two aspects of the problem
are studied. First is, as above, to find minimal value of
\eqref{shortfall risk} which will be denoted here by $\Phi_1^l(x)$
and called minimal risk function. Second is to minimize initial
costs when \eqref{shortfall risk} is smaller or equal to $v$. The
corresponding cost minimizing function is denoted by $\Phi_2^l(v)$.
Proposition 3.1 and Theorem 3.2 in \cite{FL2} provide description of 
solution to the first problem in a general semimartingale
framework and Section 7 in \cite{FL2} shows its relation to the solution 
of the second problem.  In fact these results can be treated as a general 
method for finding $\Phi_1^l(x)$, $\Phi_2^l(v)$, but they do not provide  
explicit formulas in general situation. Then using regularity of the one dimensional
Black-Scholes model both problems have been solved explicitly for a call
option, see Section 6 of \cite{FL2}.

In this paper we examine a multidimensional Black-Scholes model and
extend the results of \cite{FL2} towards more direct formulas for
the functions $\Phi_1^l$ and $\Phi_2^l$.  First we
treat the case when $l$ is linear, i.e. $l(x)=x$. Using general
results from \cite{FL2} and the fact that density of the martingale
measure is regular, we show that
\begin{gather*}
 \Phi^l_1=\Psi_1\circ\Psi^{-1}_2, \qquad \Phi^l_2=\Psi_2\circ\Psi^{-1}_1,
\end{gather*}
where $\Psi_1$, $\Psi_2$ are certain deterministic function, for 
precise formulation see Theorem \ref{tw dla
linear loss function}. This shows in particular that $\Phi_2^l$ is the inverse of
$\Phi_1^l$. We show similar results for a strictly convex loss function $l$.
As an immediate consequence of Theorem 3.2 in
\cite{FL2} we obtain the
following characterization of the risk minimizing function
\begin{gather*}
 \Phi^l_1=\Psi^l_1\circ(\Psi^{l}_2)^{-1},
\end{gather*}
where again $\Psi^l_1$, $\Psi^l_2$ are certain deterministic
functions. Analogous result for the function $\Phi^l_2$ requires
an auxiliary result - Proposition
\ref{tw auxilliary dla Psi^l_2}. Finally,
in Theorem \ref{tw dla Psi^l_2} we show that
\begin{gather*}
 \Phi^l_2=\Psi^l_2\circ(\Psi^{l}_1)^{-1}.
\end{gather*}
The risk functions $\Phi_1, \Phi_2,$ $\Phi^l_1, \Phi^l_2$ are thus determined provided
that the auxiliary functions $\Psi_1, \Psi_2$, $\Psi^l_1, \Psi^l_2$
are given. We present concrete integral forms of these functions for some widely traded
derivatives like digital option, quantos, outperformance and
spread options. The case when $l$ is a linear loss function as well as 
$l(x)=x^p, p>1$ is treated.

Let us stress the fact that both functions $\Phi^l_1,\Phi^l_2$
reflect the interplay between hedging risk and trading costs and 
thus they serve as an important tool for risk management. Although the
model under consideration is a particular case of a general
framework studied in \cite{FL2}, it is commonly used in practice due
to its tractability, see \cite{Glasserman} p.104. Thus more explicit
computing methods for finding $\Phi_1^l$ and $\Phi_2^l$ seem to be
important for practitioners. The results presented in the paper 
extend the results from \cite{Barski}, where analogous integral representations for 
the quantile hedging problem for basket derivatives have been shown.

The paper is organized as follows. In Section \ref{Problem
formulation} we describe the model settings and formulate
the problem strictly. Section \ref{Main results} contains the main results
which consist of two parts concerning a linear and a convex loss
function respectively. Section \ref{Two dimensional model} is
devoted to presenting explicit integral form for the risk functions in two dimensional
model when $l(x)=x$ and $l(x)=x^p$ with $p>1$.

\section{Problem formulation}\label{Problem formulation}
Let $(\Omega,\mathcal{F},\mathcal{F}_{t}, t\in[0,T],P)$ be a filtered probability space supporting a $d$-dimensional 
Wiener process $W=(W^1,W^2,...,W^d)$ with a positive definite correlation matrix
of the form
\begin{gather*}
Q=\left[\begin{array}{rrrrr}
1&\rho_{1,2}&\rho_{1,3}&\ldots&\rho_{1,d}\\
\rho_{2,1}&1&\rho_{2,3}&\ldots&\rho_{2,d}\\
\vdots&\vdots&\vdots&\vdots&\vdots\\
\rho_{d,1}&\rho_{d,2}&\rho_{d,3}&\ldots&1\\
\end{array}\right],
\end{gather*}
where
\begin{gather*}
\rho_{i,j}=cor\left\{W^i_1,W^j_1\right\}, \quad i,j=1,2,...,d.
\end{gather*}
The process $W$ as above will be called a $Q$-Wiener process.
The multidimensional Black-Scholes model is specified by the dynamics of $d$ stocks,
\begin{gather*}
dS^i_t=S^i_t(\alpha_i dt+\sigma_{i}dW^i_t), \qquad i=1,2,...,d, \quad
t\in[0,T],
\end{gather*}
and evolution of the money market account
\begin{gather*}
dB_t=rB_tdt, \quad t\in[0,T].
\end{gather*}
Above $\alpha_i\in\mathbb{R}$, $\sigma_i>0, \ i=1,2,...,d$ and $r$ stands 
for a constant interest rate. It is known that such a
market is complete and that the unique martingale measure
$\tilde{P}$ is given by the density
\begin{gather}\label{gestosc miary martyngalowej}
\frac{d\tilde{P}}{dP}=\tilde{Z}_T:=e^{-(Q^{-1}[\frac{\alpha-r\mathbf{1}_d}{\sigma}],W_T)-\frac{1}{2}\mid
Q^{-\frac{1}{2}}[\frac{\alpha-r\mathbf{1}_d}{\sigma}]\mid^2T},
\qquad t\in[0,T],
\end{gather}
with the notation
\begin{gather*}
Q^{-1}\left[\frac{\alpha-r\mathbf{1}_d}{\sigma}\right]:=Q^{-1}\left[\begin{array}{r}
\frac{\alpha_1-r}{\sigma_1}\\
\frac{\alpha_2-r}{\sigma_2}\\
\vdots\\
\frac{\alpha_d-r}{\sigma_d}\\
\end{array}\right], \qquad t\in[0,T],
\end{gather*}
for more details see, for instance, \cite{Barski}. Moreover,
\begin{gather*}
\widetilde{W}_t:=W_t+\frac{\alpha-r\mathbf{1}_d}{\sigma} \ t, \qquad
t\in[0,T],
\end{gather*}
is a $Q$- Wiener process under $\widetilde{P}$. The dynamics of the
prices under the measure $\widetilde{P}$ can be written as
\begin{gather*}
dS^i_t=S^i_t(rdt+\sigma_id\widetilde{W}^i_t), \quad i=1,2,...,d, \quad t\in[0,T].
\end{gather*}
The wealth process corresponding to the initial endowment $x$ and
the trading strategy $\pi$ is given by
\begin{gather*}
X^{x,\pi}_0=x, \quad
X^{x,\pi}_t:=\pi^{0}_{t}B_t+\sum_{i=1}^{d}\pi^{i}_{t}S^{i}_{t},
\quad t\in[0,T].
\end{gather*}
Each strategy is assumed to be admissible, i.e. $X^{x,\pi}_{t}\geq0$
for each $t\in[0,T]$ almost surely and self-financing, i.e.
\begin{gather*}
dX^{x,\pi}_{t}=\pi^{0}_{t}dB_t+\sum_{i=1}^{d}\pi^{i}_{t}dS^{i}_{t},
\quad t\in[0,T].
\end{gather*}
A contingent claim is represented by an $\mathcal{F}_T$- measurable
random variable $H$ which is assumed to be nonnegative, i.e.
$H\geq0$ and $\mathbf{E}[e^{-rT}\tilde{Z}_TH]<+\infty$. As the market is complete, the price of $H$ defined by
\begin{gather*}
p(H):=\inf\left\{x: \exists \pi \ \text{s.t.} \ P(X^{x,\pi}_T\geq
H)=1\right\}
\end{gather*}
is given by $p(H)=\mathbf{\tilde{E}}[e^{-rT}H]=\mathbf{E}[e^{-rT}\tilde{Z}_TH]$.

Trader's attitude towards risk is measured by
\begin{gather*}
\mathbf{E}[l((H-X_T^{x,\pi})^{+})],
\end{gather*}
where $l:[0,+\infty)\longrightarrow[0,+\infty)$ is a loss function
which is assumed to be increasing with $l(0)=0$. It is clear that if
$x\geq p(H)$ then the risk equals zero for the replicating strategy.
In the opposite case the risk is strictly positive and the question
under consideration is to find a strategy such that
\begin{gather*}
\mathbf{E}[l((H-X_T^{x,\pi})^{+})]\underset{\pi}{\longrightarrow}\min.
\end{gather*}
We will refer the corresponding function
$\Phi_1:[0,+\infty)\longrightarrow[0,\mathbf{E}[l(H)]]$ given by
\begin{gather}\label{definicja Psi_1}
\Phi^l_1(x):=\min_{\pi}\mathbf{E}[l((H-X_T^{x,\pi})^{+})],
\end{gather}
as the {\it minimal risk function}. The strategy $\hat{\pi}$ such
that $\mathbf{E}[l((H-X_T^{x,\hat{\pi}})^{+})]=\Phi^l_1(x)$ will be called the
{\it risk minimizing strategy for $x$}. If
$x\geq p(H)$ then $\Phi^l_1(x)=0$ and $\Phi^l_1(x)>0$ otherwise.

We also consider the cost reduction problem. Let $v\geq 0$ be a
fixed number describing the level of shortfall risk accepted by the
trader. We are searching for a minimal initial cost such that there
exists a strategy with the risk not exceeding $v$, i.e.
\begin{gather*}
x\longrightarrow \min; \quad \text{$\exists$ $\pi$ \ s.t.} \quad
\mathbf{E}[l((H-X_T^{x,\pi})^{+})]\leq v.
\end{gather*}
The {\it cost reduction function}
$\Phi^l_2:[0,+\infty)\longrightarrow[0,p(H)]$ is thus defined
by
\begin{gather}\label{definicja Psi_2}
\Phi^l_2(v):=\min\left\{x: \exists \pi  \ \text{s.t.} \
\mathbf{E}[l((H-X_T^{x,\pi})^{+})]\leq v \right\}.
\end{gather}
The strategy $\hat{\pi}$ such that
$\mathbf{E}[l((X^{\Phi^l_2(v),\hat{\pi}}_T-H)^{+})]\leq v$ will be
called the {\it cost minimizing strategy for $v$}. Notice that
$\Phi^l_2(0)=p(H)$.

\section{Main results}\label{Main results}
\subsection{Linear loss function}
In this section we examine the case when $l(x)=x$.  Denote for simplicity the
corresponding functions $\Phi_1^l, \Phi_2^l$ by $\Phi_1, \Phi_2$
respectively. 

Let us start with two auxiliary results.
\begin{lem}\label{lemat o funkcji g}
Let $X\geq0,Y\geq0$ be random variables such that $\mathbf{E} X<+\infty$. Then the function $g:[0,+\infty)\rightarrow [0,+\infty)$ given by
\begin{gather*}
g(c):=\mathbf{E}[X\mathbf{1}_{\{Y\geq c\}}]
\end{gather*}
\begin{enumerate}[a)]
\item is left continuous on $(0,+\infty)$ with right limits on $[0,+\infty)$,
\item is right continuous on $[0,+\infty)$ if the cumulative distribution function of $Y$ is continuous,
\item is strictly decreasing if for any $0\leq a<b<+\infty$ holds
\begin{gather}\label{warunek c}
P\left(X>0, Y\in[a,b)\right)>0.
\end{gather}
\end{enumerate}
\end{lem}
{\bf Proof:}\\
The function $g$ is decreasing and thus it has right
and left limits. If $X=0$ then the assertion follows trivially. In the opposite case
let us consider an auxiliary probability measure
$\hat{P}$ defined by
\begin{gather*}
\frac{d\hat{P}}{dP}=\frac{X}{\mathbf{E}[X]},
\end{gather*}
which is absolutely continuous with respect to $P$, i.e. $\hat{P} \ll P$.\\

\noindent
$a)$ For any $c>0$ we have
\begin{gather*}
\bigcap_{n}\left\{c-\frac{1}{n}\leq Y<c\right\}=\emptyset,
\end{gather*}
and thus
\begin{gather*}
\mid g(c-\frac{1}{n})-g(c)\mid=\mathbf{E}\left(X\mathbf{1}_{\{c-\frac{1}{n}\leq Y<c\}}\right)=\mathbf{E}[X]\hat{P}\left(c-\frac{1}{n}\leq Y<c\right)\underset{n}{\longrightarrow}0.
\end{gather*}
$b)$ For any $c\geq 0$ we have
\begin{gather*}
\bigcap_{n}\left\{c\leq Y<c+\frac{1}{n}\right\}=\left\{Y=c\right\},
\end{gather*}
and thus
\begin{align*}
\mid g(c)-g(c+\frac{1}{n})\mid&=\mathbf{E}\left(X\mathbf{1}_{\{c\leq Y<c+\frac{1}{n}\}}\right)\\[1ex]
&=\mathbf{E}[X]\hat{P}\left(c\leq Y<c+\frac{1}{n}\right)\underset{n}{\longrightarrow}\mathbf{E}[X]\hat{P}(Y=c)=0,
\end{align*}
as $\hat{P} \ll P$ and $P(Y=c)=0$.\\
$c)$ Let us notice that \eqref{warunek c} is equivalent to the condition
\begin{gather*}
\exists \varepsilon>0 \ \text{s.t.} \ P(X>\varepsilon, Y\in[a,b))>0,
\end{gather*}
and thus for $0\leq a<b<+\infty$ we have
\begin{align*}
\mid g(a)-g(b)\mid&=\mathbf{E}\left(X\mathbf{1}_{\{a \leq Y<b\}}\right)\\
&=\mathbf{E}\left(X\mathbf{1}_{\{a \leq Y<b\}}\mathbf{1}_{\{X=0\}}\right)+
\mathbf{E}\left(X\mathbf{1}_{\{a \leq Y<b\}}\mathbf{1}_{\{X>0\}}\right)\\
&\geq\mathbf{E}\left(X\mathbf{1}_{\{a \leq Y<b\}}\mathbf{1}_{\{X>\varepsilon\}}\right)\geq\varepsilon P(X>\varepsilon, a\leq Y<b)>0.
\end{align*}
\hfill$\square$

\begin{rem}
Let us notice that condition \eqref{warunek c} implies that $Y$ has a strictly increasing cumulative distribution function. 
\end{rem}

\begin{cor}\label{corollary o funkcji g}
If the cumulative distribution function of $Y$ is continuous then the function $g$ in Lemma \ref{lemat o funkcji g} is continuous on $(0,+\infty)$ and right continuous at $0$.
\end{cor}

\begin{prop}\label{prop o wektorze normalnym}
Let $(Z_1,Z_2)$ be a random vector with nondegenerate normal distribution on a plane. Let
$f,g$ be functions such that
\begin{gather*}
 f:\mathbb{R}^2\longrightarrow(0,+\infty),\\[1ex]
 h:\mathbb{R}\longrightarrow(0,+\infty) \quad \text{is strictly monotone}.
\end{gather*}
Let $\alpha,\beta,\gamma,\delta\in\mathbb{R}$ be such that the vectors $(\alpha,\beta),(\gamma,\delta)$ are not parallel, i.e. $(\alpha,\beta)\nparallel(\gamma,\delta)$. Let
\begin{align*}
 X:=f(Z_1,Z_2)\mathbf{1}_{\{\alpha Z_1+\beta Z_2>k\}}, \quad Y:=h(\gamma Z_1+\delta Z_2),
\end{align*}
where $k$ is some constant. Then the function $g(c):=\mathbf{E}[X\mathbf{1}_{\{Y\geq c\}}]$ is 
strictly decreasing on $[0,+\infty)$.
\end{prop}
{\bf Proof:} We will show that \eqref{warunek c} in Lemma \ref{lemat o funkcji g} holds. We have
\begin{gather*}
P(X>0,Y\in[a,b))=P\Big(\alpha Z_1+\beta Z_2>k, h^{-1}(a)\leq\gamma Z_1+\delta Z_2<h^{-1}(b)\Big)
\end{gather*}
for the case when $h$ is strictly increasing. The probability above
is positive because the set
\begin{gather*}
 \Big\{(x,y): \alpha x+\beta y>k, h^{-1}(a)\leq\gamma x+\delta y<h^{-1}(b)\Big\}
\end{gather*}
is of positive Lebesgue measure and $(Z_1,Z_2)$ has a nondegenerate distribution.
\hfill$\square$

\vskip2ex

In the sequel we will need  two auxiliary functions defined by
\begin{align}\label{Psi_1}
\Psi_1(c)&:=\mathbf{E}(H\mathbf{1}_{A_c}),\\[2ex]
\Psi_2(c)&:=\mathbf{\tilde{E}}(H\mathbf{1}_{A_c}),\label{Psi_2}
\end{align}
where
\begin{gather*}
A_c:=\{\tilde{Z}^{-1}_T\geq c\}, \quad c\geq 0,
\end{gather*}
and $\tilde{Z}_T$ is given by \eqref{gestosc miary martyngalowej}.
Let us notice that due to the fact that $Q$ is nonsingular the random variable
\begin{gather}
\tilde{Z}^{-1}_T:=e^{(Q^{-1}[\frac{\alpha-r\mathbf{1}_d}{\sigma}],W_T)+\frac{1}{2}\mid
Q^{-\frac{1}{2}}[\frac{\alpha-r\mathbf{1}_d}{\sigma}]\mid^2T}
\end{gather}
has a continuous cumulative distribution function with respect to $P$ and $\tilde{P}$. Thus it follows from Corollary \ref{corollary o funkcji g} that the functions $\Psi_1$, $\Psi_2$ are continuous for any $H\geq0$. It is clear that both are decreasing. For some special contingent claims the functions are strictly decreasing.
Indeed, using Proposition \ref{prop o wektorze normalnym} one can show that this is the case for the following payoffs.
\begin{ex}
The functions $\Psi_1$, $\Psi_2$ are strictly decreasing if
\begin{enumerate}[a)] 
 \item $H$ is a digital option, i.e. $H=K\mathbf{1}_{\{S^1_T\geq S^2_T\}}$ and $(\sigma_1,-\sigma_2)\nparallel Q^{-1}[\frac{\alpha-r\mathbf{1}_d}{\sigma}]$,
 \item $H$ is a quanto domestic option, i.e. $H=S^2_T(S^1_T-K)^{+}$ and $(\sigma_1,0)\nparallel Q^{-1}[\frac{\alpha-r\mathbf{1}_d}{\sigma}]$,
 \item $H$ is a quanto foreign option, i.e. $H=(S^1_T-\frac{K}{S^2_T})^{+}$ and $(\sigma_1,\sigma_2)\nparallel Q^{-1}[\frac{\alpha-r\mathbf{1}_d}{\sigma}]$.
\end{enumerate}
\end{ex}

\noindent
Below we present the description of the risk functions $\Phi_1$, $\Phi_2$.

\begin{tw}\label{tw dla linear loss function}
\begin{enumerate}[a)]
\item Let $c=c(x)$ be a solution of the equation
\begin{gather}\label{istnienie stalej dla Psi_2}
\Psi_2(c)=e^{rT}x, \qquad x\in[0,p(H)).
\end{gather}
Then
\begin{gather*}
\Phi_1(x)=
\begin{cases}
\Psi_1(0)-\Psi_1(c) \ &\text{for}\ x\in[0,p(H)),\\
0 \ &\text{for} \  x\geq p(H).
\end{cases}
\end{gather*}
Moreover, the replicating strategy for the payoff $H\mathbf{1}_{A_{c(x)}}$ is a risk minimizing strategy for $x$.

\item Let $c=c(v)$ be a solution of the equation
\begin{gather}\label{istnienie stalej dla Psi_1}
\Psi_1(c)=\Psi_1(0)-v, \qquad v\in[0,\mathbf{E}[H]).
\end{gather}
Then
\begin{gather*}
\Phi_2(v)=
\begin{cases}
e^{-rT}\Psi_2(c) \ &\text{for} \  v\in[0,\mathbf{E}[H]),\\
0 \ &\text{for}\ v\geq\mathbf{E}[H].
\end{cases}
\end{gather*}
Moreover, the replicating strategy for the payoff
$H\mathbf{1}_{A_{c(v)}}$ is a cost minimizing strategy for $v$.
\end{enumerate}
\end{tw}
{\bf Proof:}
First let us notice that the equations \eqref{istnienie stalej dla Psi_2}, \eqref{istnienie stalej dla Psi_1} actually have solutions.
Indeed, it follows from the fact that $\Psi_1,\Psi_2$ are continuous and decreasing with images
$[0,\mathbf{E}[H]]$, $[0,e^{rT}p(H)]$ respectively.

For any admissible strategy $(x,\pi)$ let us define the success function
\begin{gather*}
 \varphi_{x,\pi}:=\mathbf{1}_{\{X^{x,\pi}_{T}\geq H\}}+\frac{X^{x,\pi}_{T}}{H}\mathbf{1}_{\{X^{x,\pi}_{T}< H\}}.
\end{gather*}
One can check the following identity
\begin{gather*}
(H-X^{x,\pi}_{T})^{+}=H-X^{x,\pi}_{T}\wedge H=H-H\varphi_{x,\pi},
\end{gather*}
which implies that
\begin{gather}\label{wzor na expected loss w terminach varphi}
\mathbf{E}[(H-X^{x,\pi}_{T})^{+}]=\mathbf{E}[H]-\mathbf{E}[H\varphi_{x,\pi}].
\end{gather}

\noindent
$a)$ In view of \eqref{wzor na expected loss w terminach varphi} the problem \eqref{definicja Psi_1} of finding $\Phi_1(x)$ is equivalent to that of finding the strategy
$\pi$ satisfying
\begin{align*}
\mathbf{E}[H\varphi_{x,\pi}]\underset{\pi}{\rightarrow} \max.
\end{align*}
If $x\geq p(H)$ then $\varphi_{x,\pi}=1$ for the replicating
strategy and $\Phi_1(x)=0$, so consider the case $0\leq x<p(H)$. Let
us formulate an auxiliary problem of determining
$\varphi\in\mathcal{R}$ solving
\begin{gather}\label{auxiliary problem 1}
\begin{cases}
\mathbf{E}[H\varphi]{\rightarrow} \max,\\[2ex]
\mathbf{\tilde{E}}[e^{-rT}H\varphi]\leq x,
\end{cases}
\end{gather}
where
\begin{gather}\label{definicja R}
\mathcal{R}:=\{\varphi: 0\leq\varphi\leq1 \ \text{and} \ \varphi \ \text{is} \ \mathcal{F}_T-\text{measurable}\}.
\end{gather}
It is clear that if $\hat{\varphi}$ such that $\mathbf{\tilde{E}}[e^{-rT}H\hat{\varphi}]=x$ is a solution of \eqref{auxiliary problem 1} then
the replicating strategy $\tilde{\pi}$ for the payoff $H\hat{\varphi}$ is a risk minimizing strategy for $x$ and
\begin{gather}\label{wzor na Psi_1 przy uzyciu varphi}
\Phi_1(x)=\mathbf{E}[(H-X_{T}^{x,\tilde{\pi}})^{+}]=\mathbf{E}[H]-\mathbf{E}[H\hat{\varphi}].
\end{gather}
Thus now let us focus on determining solution $\hat{\varphi}$ of
\eqref{auxiliary problem 1}. To this end introduce two probability
measures $P_1,P_2$ with densities
\begin{gather*}
\frac{dP_1}{dP}=\frac{H}{\mathbf{E}[H]}, \quad \frac{dP_2}{dP}=\frac{e^{-rT}\tilde{Z}_TH}{\mathbf{E}[e^{-rT}\tilde{Z}_TH]} .
\end{gather*}
Then \eqref{auxiliary problem 1} reads as
\begin{gather}\label{auxiliary problem 2}
\begin{cases}
\mathbf{E}^{P_1}[\varphi]{\longrightarrow} \max,\\[2ex]
\mathbf{E}^{P_2}[\varphi]\leq \frac{x}{p(H)},
\end{cases}
\end{gather}
which is a standard problem in the theory of statistical tests. One should try to search for
the solution in the class of $0-1$ valued functions of the form $\mathbf{1}_{A_c}; c\geq0$, where
\begin{align*}
A_c&:=\Big\{\frac{dP_1}{dP_2}\geq c\Big\}=\Big\{\frac{dP_1}{dP}\frac{dP}{dP_2}\geq c\Big\}=\Big\{\frac{H}{\mathbf{E}[H]}\frac{\mathbf{E}[\tilde{Z}_TH]}{\tilde{Z}_TH}\geq c\Big\}\\[1ex]
&=\Big\{\tilde{Z}_T^{-1}\geq c\frac{\mathbf{E}[H]}{\mathbf{E}[\tilde{Z}_TH]}\Big\}.
\end{align*}
For the sake of simplicity we can reparametrize $A_c$ by denoting the constant $c\frac{\mathbf{E}[H]}{\mathbf{E}[\tilde{Z}_TH]}$ above just by $c$. Then $A_c$ is of the form
\begin{gather*}
A_c:=\Big\{\tilde{Z}_T^{-1}\geq c\Big\}.
\end{gather*}
It is known by the Neyman-Pearson lemma that if there exists $c=c(x)$ such that
\begin{gather}\label{warunek na c 1}
\mathbf{E}^{P_2}[\mathbf{1}_{A_{c}}]=P_2(A_{c})=\frac{x}{p(H)},
\end{gather}
then the solution of \eqref{auxiliary problem 2}, or equivalently \eqref{auxiliary problem 1}, is given by $\hat{\varphi}=\mathbf{1}_{A_{c(x)}}$. But let us notice that \eqref{warunek na c 1} is equivalent to the following
\begin{gather*}
\Psi_2(c)=e^{rT}x,
\end{gather*}
and the existence of the required constant $c$ follows from
\eqref{istnienie stalej dla Psi_2}. Finally, coming back to
\eqref{wzor na Psi_1 przy uzyciu varphi} and using definition of
$\Psi_1$, we  obtain
\begin{gather*}
\Phi_1(x)=\mathbf{E}[H]-\mathbf{E}[H\hat{\varphi}]=\mathbf{E}[H]-\mathbf{E}[H\mathbf{1}_{A_c}]=\Psi_1(0)-\Psi_1(c).
\end{gather*}

\noindent $b)$ If $v\geq\mathbf{E}[H]$ then the cost minimizing
strategy is trivial, i.e. $(x=0,\pi=0)$ and thus $\Phi_2(v)=0$. Let
us focus on the case when $v\in[0,\mathbf{E}[H])$. In view of
\eqref{wzor na expected loss w terminach varphi} the cost minimizing
strategy is the one which solves the problem
\begin{gather*}
 \begin{cases}
  \mathbf{E}[H\varphi_{x,\pi}]\geq\mathbf{E}[H]-v\\[2ex]
    \mathbf{\tilde{E}}[e^{-rT}H\varphi_{x,\pi}]\longrightarrow\min.
 \end{cases}
\end{gather*}
We are thus looking for a solution $\hat{\varphi}\in\mathcal{R}$ of the problem
\begin{gather}\label{auxiliary problem 3}
 \begin{cases}
  \mathbf{E}[H\varphi]\geq\mathbf{E}[H]-v\\[2ex]
    \mathbf{\tilde{E}}[e^{-rT}H\varphi]\longrightarrow\min.
 \end{cases}
\end{gather}
If \eqref{auxiliary problem 3} has a solution satisfying
$\mathbf{E}[H\hat{\varphi}]=\mathbf{E}[H]-v$ then the cost
minimizing strategy is that one which replicates $H\hat{\varphi}$ and
the cost minimizing function equals
\begin{gather}\label{wzor na Phi_2 przy uzyciu varphi}
 \Phi_2(r)=e^{-rT}\mathbf{\tilde{E}}[H\hat{\varphi}].
\end{gather}
Let us focus on determining the solution $\hat{\varphi}$ of \eqref{auxiliary problem 3}. Using notation from the part $(a)$ we can reformulate \eqref{auxiliary problem 3} to the form
\begin{gather}\label{auxiliary problem 4}
\begin{cases}
 \mathbf{E}^{P_1}[\varphi]\geq\frac{\mathbf{E}[H]-v}{\mathbf{E}[H]}\\[2ex]
\mathbf{E}^{P_2}[\varphi]\longrightarrow \min.
\end{cases}
 \end{gather}
It can be shown in the same way as in the proof of Neyman-Pearson lemma that the solution should be searched in the
$0-1$ valued functions of the form $\mathbf{1}_{B_c}; c\geq0$, where
\begin{gather*}
 B_c:=\Big\{\frac{dP_2}{dP_1}\leq c\Big\}=\Big\{\frac{dP_2}{dP}\frac{dP}{dP_1}\leq c\Big\}=\Big\{\tilde{Z}^{-1}_T\geq\frac{1}{c}\frac{\mathbf{E}[H]}{\mathbf{E}[\tilde{Z}_TH]}\Big\}.
\end{gather*}
Denoting, for simplicity, the constant $\frac{1}{c}\frac{\mathbf{E}[H]}{\mathbf{E}[\tilde{Z}_TH]}$ above by $c$, we have
\begin{gather*}
 B_c=\{\tilde{Z}_T^{-1}\geq c\}.
\end{gather*}
If there exists constant $c=c(v)$ satisfying
\begin{gather}\label{warunek na c(r)}
\mathbf{E}^{P_1}[\mathbf{1}_{B_c}]=P_1(B_c)=\frac{\mathbf{E}[H]-v}{\mathbf{E}[H]}
\end{gather}
then $\hat{\varphi}=\mathbf{1}_{B_c}$ is a solution of \eqref{auxiliary problem 4} or, equivalently, \eqref{auxiliary problem 3}.
Let us notice that \eqref{warunek na c(r)} can be written as
\begin{gather*}
\Psi_1(c)=\Psi_1(0)-v
\end{gather*}
and existence of the required constant $c(v)$ follows from
\eqref{istnienie stalej dla Psi_1}. Coming back to \eqref{wzor na
Phi_2 przy uzyciu varphi} we obtain
\begin{gather*}
\Phi_2(v)=e^{-rT}\mathbf{\tilde{E}}[H\mathbf{1}_{B_c}]=e^{-rT}\Psi_2(c).
\end{gather*}
\hfill $\square$

\subsection{Convex loss function}
In this section we study the case when
$l:[0,+\infty)\longrightarrow[0,+\infty)$ is an increasing, strictly
convex function such that $l(0)=0$. We assume that $l\in
C^{2}(0,+\infty)$ and that $l^{\prime}$ is strictly increasing with
$l^{\prime}(0+)=0$, $l^{\prime}(+\infty)=+\infty$. The inverse of
the first derivative will be denoted by $I$, i.e.
\begin{gather*}
I=(l^{\prime})^{-1}.
\end{gather*}
Moreover, the contingent claim $H$ is assumed to satisfy $\mathbf{E}[l(H)]<+\infty$.
The functions $\Phi^l_1, \Phi^l_2$ can be characterized in terms of
the functions
\begin{align}\label{Psi_1^l}
\Psi^l_1(c)&:=\mathbf{E}[l((1-\varphi_c)H)]\\[2ex]\label{Psi_2^l}
\Psi^l_2(c)&:=\mathbf{\tilde{E}}[H\varphi_c],
\end{align}
where $\varphi_c$ is defined by
\begin{gather}\label{postac varphi_c}
\varphi_c:=\left\{1-\left(\frac{I(c\tilde{Z}_T)}{H}\wedge
1\right)\right\}\mathbf{1}_{\{H>0\}}, \qquad c\geq 0.
\end{gather}
It was shown in \cite{FL2}, Theorem 5.1, that the problem of
determining $\Phi^l_1$ is equivalent to finding the solution $\tilde{\varphi}$
of the problem
\begin{gather}\label{auxiliary problem 1l}
\begin{cases}
\mathbf{E}[l((1-\varphi)H)]\underset{\varphi\in\mathcal{R}}{\longrightarrow} \min\\[2ex]
\mathbf{\tilde{E}}[e^{-rT}H\varphi]\leq x,
\end{cases}
\end{gather}
where $\mathcal{R}$ is defined in \eqref{definicja R}. Then
$\Phi^l_1(x)=\mathbf{E}[l((1-\tilde{\varphi})H)]$ and the risk
minimizing strategy is the one which replicates $H\tilde{\varphi}$.
Moreover, since the function $\Psi^l_2$ is continuous with the image
$[0,e^{rT}p(H)]$, see the proof of Theorem 5.1 in \cite{FL2}, it
follows that for any $x\in[0,e^{rT}p(H)]$ there exists constant $c$
such that $\Psi^l_2(c)=\mathbf{\tilde{E}}[H\varphi_c]=e^{rT}x$. Such
$\varphi_c$ solves the auxiliary problem \eqref{auxiliary problem
1l} and thus
\begin{gather*} \Phi^l_1(x)=\mathbf{E}[l((1-\varphi_c)H)],
\end{gather*}
and the minimal risk strategy is that replicating the payoff
$H\varphi_c$, see Theorem 3.2 in \cite{FL2}. Thus the results from
\cite{FL2} can be expressed in our notation as follows.
\begin{tw}\label{tw dla Psi^l_1}
Let $c=c(x)$ be a solution of the equation
\begin{gather*}
\Psi^l_2(c)=e^{rT}x, \quad x\in[0,p(H)).
\end{gather*}
Then
\begin{gather*}
\Phi^l_1(x)=
\begin{cases}
\Psi^l_1(c) \ &\text{for}\ x\in[0,p(H)),\\
0 \ &\text{for} \  x\geq p(H).
\end{cases}
\end{gather*}
\end{tw}

Although Theorem \ref{tw dla Psi^l_1} is only a reformulation of
Theorem 3.2 in \cite{FL2}, it provides an effective method for
practical applications if one is able to derive the functions
$\Psi^l_1$, $\Psi^l_2$ for concrete derivatives.

In the sequel we will show that the function $\Phi^l_2$ can be
characterized in terms of the functions $\Psi^l_1, \Psi^l_2$ as
well. It is easy to show that the cost reduction problem is
equivalent to that of finding $\varphi\in\mathcal{R}$ such that
\begin{gather}\label{auxiliary problem 2l}
\begin{cases}
\mathbf{E}[l((1-\varphi)H)]\leq v\\[2ex]
\mathbf{\tilde{E}}[e^{-rT}H\varphi]\longrightarrow \min.
\end{cases}
\end{gather}
Let us notice that \eqref{auxiliary problem 2l} can not be solved
with the same method as \eqref{auxiliary problem 1l}. In
\eqref{auxiliary problem 1l} the constraints are linear and thus the
solution could be found via Neyman-Pearson approach to the
variational problem, see the proof of Theorem 5.1 in \cite{FL2} and
p.210 in \cite{Karlin}. The constraints in \eqref{auxiliary problem 2l}
are no longer linear and the method above fails. Below we present the
proof based on Lagrange multipliers.

\begin{prop}\label{tw auxilliary dla Psi^l_2}
Let $l^{\prime\prime}$ be increasing and let $H$ additionally satisfy\\
$\mathbf{E}[l^{\prime}(H)H]<+\infty$ and $\mathbf{E}[l^{\prime\prime}(H)H^2]<+\infty$.
Then the solution of the problem \eqref{auxiliary problem 2l} is of the form
\begin{gather*}
\tilde{\varphi}:=\left\{1-\left(\frac{I(c\tilde{Z}_T)}{H}\wedge 1\right)\right\}\mathbf{1}_{\{H>0\}}
\end{gather*}
where $c$ is such that $\mathbf{E}[l((1-\tilde{\varphi})H)]=v$.
\end{prop}
{\bf Proof:} First let us notice that if $\varphi\in\mathcal{R}$ is
a solution to \eqref{auxiliary problem 2l} then necessarily
$\mathbf{E}[l((1-\varphi)H)]=v$. Indeed, assume to the contrary that
$\varphi$ is a solution to \eqref{auxiliary problem 2l} with
$\mathbf{E}[l((1-\varphi)H)]<v$ and consider a family of random
variables $\varphi_{\alpha}:=\varphi\wedge\alpha; \alpha\in[0,1]$.
Then the function $\alpha\rightarrow
\mathbf{E}[l((1-\varphi_{\alpha})H)]$ is continuously decreasing
from $\mathbf{E}[l(H)]$ to $0$. Thus there exists
$\tilde{\alpha}\in[0,1]$ such that
$\mathbf{E}[l((1-\varphi_{\tilde{\alpha}})H)]=v$. Then
$\varphi_{\tilde{\alpha}}\leq\varphi$ and thus
$\mathbf{\tilde{E}}[H\varphi_{\tilde{\alpha}}]<\mathbf{\tilde{E}}[H\varphi]$,
which is a contradiction.

Let $\varphi\neq\tilde{\varphi}$ be any element of $\mathcal{R}$
such that $\mathbf{E}[l((1-\varphi)H)]=v$. We need to show that
$\mathbf{\tilde{E}}[H\tilde{\varphi}]\leq\mathbf{\tilde{E}}[H\varphi]$.
Let us define $\varphi_{\varepsilon}$ by
\begin{gather*}
\varphi_{\varepsilon}:=(1-\varepsilon)\tilde{\varphi}+\varepsilon\varphi,\qquad \varepsilon\in[0,1],
\end{gather*}
and the function
\begin{gather*}
  F_{\varphi}(\varepsilon):=\mathbf{\tilde{E}}(H\varphi_{\varepsilon})=\mathbf{E}(\tilde{Z}_TH\varphi_{\varepsilon}).
\end{gather*}
We need to show that $F_{\varphi}(0)\leq F_{\varphi}(1)$. We will
show that $F_{\varphi}$ has minimum at $0$. Let us define the
auxiliary function
\begin{align*}
 G_{\varphi}(\varepsilon):=\mathbf{E}[l((1-\varphi_{\varepsilon})H)],
\end{align*}
and notice that due to the convexity of $l$ we have
$G_{\varphi}(\varepsilon)\leq v$ for each $\varepsilon\in[0,1]$.
Thus the problem of minimizing $F_{\varepsilon}$ on $[0,1]$ is
equivalent to the following
\begin{align}\label{problem F,G}
\begin{cases}
&F_{\varphi}(\varepsilon)\longrightarrow \min\\[1ex]
&G_{\varphi}(\varepsilon)\leq v,\\[1ex]
&\varepsilon\geq 0,\\[1ex]
&1-\varepsilon\geq0.
\end{cases}
\end{align}

\noindent
In view of the assumptions on $l$ and $H$, both functions $F_{\varphi}, G_{\varphi}$ are smooth with
\begin{align*}
 F^{\prime}_{\varphi}(\varepsilon)&\equiv \mathbf{E}[\tilde{Z}_T(\varphi-\tilde{\varphi})H],\\[1ex]
 G^{\prime}_{\varphi}(\varepsilon)&=\mathbf{E}[l^{\prime}((1-\varphi_{\varepsilon})H)\cdot(\tilde{\varphi}-\varphi)H],\\[1ex]
G^{\prime\prime}_{\varphi}(\varepsilon)&=\mathbf{E}[l^{\prime\prime}((1-\varphi_{\varepsilon})H)\cdot(\tilde{\varphi}-\varphi)^2H^2],
\end{align*}
and thus the Lagrange function for \eqref{problem F,G} is of the form
\begin{gather*}
L(\varepsilon,\lambda_1,\lambda_2,\lambda_3)=F_\varphi(\varepsilon)-\lambda_1(v-G_{\varphi}(\varepsilon))-\lambda_2\varepsilon-\lambda_3(1-\varepsilon).
\end{gather*}
As the function $F_{\varphi}$ is linear, it attains its minimal value at $0$ or $1$. We will show that the first and the second order differential conditions are satisfied for $\varepsilon=0$.

\noindent
The first order conditions are
{\small
\begin{gather}\label{first order conditions 1}
L^{\prime}_{\varepsilon}(\varepsilon,\lambda_1,\lambda_2,\lambda_3)=\mathbf{E}[\tilde{Z}_T(\varphi-\tilde{\varphi})H]+\lambda_1 \mathbf{E}[l^{\prime}((1-\varphi_{\varepsilon})H)\cdot(\tilde{\varphi}-\varphi)H]-\lambda_2+\lambda_3=0\\[2ex]\label{first order conditions 2}
\lambda_1,\lambda_2,\lambda_3\geq 0,  \quad
\lambda_1(v-G_{\varphi}(\varepsilon))=0, \quad
\lambda_2\varepsilon=0, \quad \lambda_3(1-\varepsilon)=0.
\end{gather}
}
By the definition of $\tilde{\varphi}$ we have
\begin{gather*}
 \tilde{\varphi}=1-\frac{I(c\tilde{Z}_T)}{H} \quad \text{and} \quad c\tilde{Z}_T=l^{\prime}((1-\tilde{\varphi})H) \quad \text{on} \ A\\[2ex]
 \tilde{\varphi}=0 \quad \text{on} \quad A^{c},
\end{gather*}
where $A:=\{c\tilde{Z}_T<l^{\prime}(H)\}$ and $A^c$ stands for the compliment of $A$. For $\varepsilon=0$ it follows from \eqref{first order conditions 2} that $\lambda_3=0$ and the equation \eqref{first order conditions 1} is of the form
\begin{align}\nonumber
&\mathbf{E}[\tilde{Z}_T(\varphi-\tilde{\varphi})H\mathbf{1}_A]+\mathbf{E}[\tilde{Z}_T(\varphi-\tilde{\varphi})H\mathbf{1}_{A^{c}}]+c\lambda_1\mathbf{E}[\tilde{Z}_T(\tilde{\varphi}-\varphi)H\mathbf{1}_A]\\[1ex]\nonumber
&+\lambda_1\mathbf{E}[l^{\prime}((1-\tilde{\varphi})H)(\tilde{\varphi}-\varphi)H\mathbf{1}_{A^c}\\[1ex] \label{warunek na lambda2}
&=(1-c\lambda_1)\mathbf{E}[\tilde{Z}_T(\varphi-\tilde{\varphi})H\mathbf{1}_{A}]+\mathbf{E}[\tilde{Z}_T\varphi
H\mathbf{1}_{A^c}]- \lambda_1\mathbf{E}[l^{\prime}(H)\varphi
H\mathbf{1}_{A^c}]=\lambda_2.
\end{align}
The left side of \eqref{warunek na lambda2} satisfies the following estimation
\begin{align*}
 &(1-c\lambda_1)\mathbf{E}[\tilde{Z}_T(\varphi-\tilde{\varphi})H\mathbf{1}_{A}]+\mathbf{E}[\tilde{Z}_T\varphi H\mathbf{1}_{A^c}]-
\lambda_1\mathbf{E}[l^{\prime}(H)\varphi H\mathbf{1}_{A^c}]\\[1ex]
&\geq (1-c\lambda_1)\mathbf{E}[\tilde{Z}_T(\varphi-\tilde{\varphi})H\mathbf{1}_{A}]+\mathbf{E}[\tilde{Z}_T\varphi H\mathbf{1}_{A^c}]-
\lambda_1 c\mathbf{E}[\tilde{Z}_T\varphi H\mathbf{1}_{A^c}]\\[1ex]
&\geq (1-c\lambda_1)\mathbf{E}[\tilde{Z}_T(\varphi-\tilde{\varphi})H\mathbf{1}_{A}+\tilde{Z}_T\varphi H\mathbf{1}_{A^c}].
\end{align*}
If $\mathbf{E}[\tilde{Z}_T(\varphi-\tilde{\varphi})H\mathbf{1}_{A}+\tilde{Z}_T\varphi H\mathbf{1}_{A^c}]>0$ then we take $\lambda_1$ such that $(1-c\lambda_1)>0$, in the opposite case, such that $(1-c\lambda_1)<0$. In both cases $\lambda_2$ given by \eqref{warunek na lambda2} is nonnegative.

The second order condition for $\varepsilon=0$ is
\begin{gather*}
 L^{\prime\prime}_{\varepsilon}(\varepsilon,\lambda_1,\lambda_2,\lambda_3)=\lambda_1\mathbf{E}[l^{\prime\prime}((1-\tilde{\varphi})H)\cdot(\tilde{\varphi}-\varphi)^2H^2]\geq0,
\end{gather*}
and thus the solution of \eqref{problem F,G} is $\varepsilon=0$.

\hfill$\square$

\noindent Proposition \ref{tw auxilliary dla Psi^l_2} and the
definitions of $\Psi^l_1, \Psi^l_2$ lead us to the following result.
\begin{tw}\label{tw dla Psi^l_2}
Assume that $l^{\prime\prime}$ is increasing and $H$ satisfies
$\mathbf{E}[l^{\prime}(H)H]<+\infty$ and $\mathbf{E}[l^{\prime\prime}(H)H^2]<+\infty$.
Let $c=c(v)$ be a solution of the equation
\begin{gather*}
 \Psi^l_1(c)=v, \quad v\in[0,\mathbf{E}[l(H)]).
\end{gather*}
Then
\begin{gather*}
\Phi^l_2(v)=
\begin{cases}
e^{-rT}\Psi^l_2(c) \ &\text{for}\ v\in[0,\mathbf{E}[l(H)]),\\
0 \ &\text{for} \  v\geq \mathbf{E}[l(H)].
\end{cases}
\end{gather*}
\end{tw}

\section{Two dimensional model}\label{Two dimensional model}
In this section we determine explicit integral formulas for the functions
$\Psi^l_1$, $\Psi^l_2$ for several popular options in the case $d=2$. Some of the results can be generalized
to higher dimensions.

 At the beginning let us introduce some notation concerned with the multidimensional normal distribution.
The fact that an $\mathbb{R}^d$-valued random vector $X$ has a normal distribution with mean $m\in\mathbb{R}^d$ and
covariance matrix $\Sigma$ will be denoted by $X\sim N_d(m,\Sigma)$ or $\mathcal{L}(X)=N_d(m,\Sigma)$.
$f_X$ stands for its density. If $d=1$ then the subscript is omitted and $N(m,\sigma)$ denotes the
normal distribution with mean $m$ and variance $\sigma$. If $X\sim
N_d(m,\Sigma)$ and $A$ is a $k\times d$ matrix then,
\begin{gather}\label{rozklad AX}
AX\sim N_k(Am,A\Sigma A^T);
\end{gather}
in particular if $a\in\mathbb{R}^{d}$ then
\begin{gather}\label{rozklad aX}
a^{T}X\sim N(a^Tm,a^T\Sigma a).
\end{gather}
Let $X$ be a random vector taking values in $\mathbb{R}^d$ and fix
an integer $0<k<d$. Let us divide $X$ into two vectors $X^{(1)}$ and
$X^{(2)}$ of lengths $k$, $d-k$ respectively, i.e.
\begin{gather*}
X^{(1)}=(X_1,X_2,...,X_k)^T, \qquad
X^{(2)}=(X_{k+1},X_{k+2},...,X_d)^T.
\end{gather*}
Analogously, divide the mean vector $m$ and the covariance matrix
$\Sigma$
\begin{displaymath}
m= \left(
\begin{array}{c}
m^{(1)}\\[1ex]
m^{(2)}
\end{array}
\right); \qquad \Sigma= \left[
\begin{array}{cc}
\Sigma^{(11)}&\Sigma^{(12)}\\[1ex]
\Sigma^{(21)}&\Sigma^{(22)}
\end{array}
\right],
\end{displaymath}
so that $\mathbf{E}X^{(1)}=m^{(1)}$, $\mathbf{E}X^{(2)}=m^{(2)}$,
$Cov X^{(1)}=\Sigma^{(11)}$, $Cov X^{(2)}=\Sigma^{(22)}$, $Cov
(X^{(1)},X^{(2)})=\Sigma^{(12)}={\Sigma^{(21)}}^T$. Denote by
$\mathcal{L}\left(X^{(1)}\mid X^{(2)}=x^{(2)}\right)$ the
conditional distribution of $X^{(1)}$ given
$X^{(2)}=x^{(2)}\in\mathbb{R}^{d-k}$. If $\Sigma^{(22)}$ is
nonsingular then
\begin{gather}\label{tw o rozkl nor warunkowym}
\mathcal{L}\left(X^{(1)}\mid
X^{(2)}=x^{(2)}\right)=N_{k}(m^{(1)}(x^{(2)}),\Sigma^{(11)}(x^{(2)})),
\end{gather}
where
\begin{align}\label{wzor na srednia i wariancje waunkowa ogolny}\nonumber
m^{(1)}(x^{(2)})&=
m^{(1)}+\Sigma^{(12)}{\Sigma^{(22)}}^{-1}(x^{(2)}-m^{(2)}), \\[2ex]
\Sigma^{(11)}(x^{(2)})&=\Sigma^{(11)}-\Sigma^{(12)}{\Sigma^{(22)}}^{-1}\Sigma^{(21)}.
\end{align}
Actually the conditional variance $\Sigma^{(11)}(x^{(2)})$ does not
depend on $x^{(2)}$ but we keep the notation for the sake of
consistency. The conditional density will be denoted by
$f_{X^{(1)}\mid X^{(2)}=x^{(2)}}(x^{(1)})$, where
$x^{(1)}\in\mathbb{R}^k$. In particular if $(X,Y)$ is a two
dimensional normal vector with parameters
\begin{displaymath}
m= \left(
\begin{array}{c}
m_1\\[1ex]
m_2
\end{array}
\right); \qquad \Sigma= \left[
\begin{array}{cc}
\sigma_{11}&\sigma_{12}\\[1ex]
\sigma_{21}&\sigma_{22}
\end{array}
\right],
\end{displaymath}
then
\begin{gather*}
\mathcal{L}(X\mid Y=y)=N(m_1(y),\sigma_1(y)),
\end{gather*}
where
\begin{gather}\label{wzor na sredia i wariancja warunkowa szczegolny}
m_1(y):= m_1+\frac{\sigma_{12}}{\sigma_{22}}(y-m_2), \quad
\sigma_1(y):= \sigma_{11}-\frac{\sigma_{12}^2}{\sigma_{22}}.
\end{gather}
If $X$ is a random vector then its distribution with respect to the measure
$\widetilde{P}$ will be denoted by $\tilde{\mathcal{L}}(X)$ and its
density by $\tilde{f}_X$. Analogously, $\tilde{f}_{X^{(1)}\mid
X^{(2)}=x^{(2)}}(x^{(1)})$ stands for the conditional density with
respect to the measure $\widetilde{P}$.

In the case case $d=2$ the correlation matrix is of the form
\begin{gather*}
Q=\left[\begin{array}{rr}
1&\rho\\
\rho&1\\
\end{array}\right],
\end{gather*}
and thus we have
\begin{gather*}
Q^{-1}=\frac{1}{\rho^2-1}\left[\begin{array}{rr}
-1&\rho\\
\rho&-1\\
\end{array}\right],\qquad
Q^{-\frac{1}{2}}=\frac{1}{2}\left[\begin{array}{rr}
\frac{1}{\sqrt{1+\rho}}+\frac{1}{\sqrt{1-\rho}}&\frac{1}{\sqrt{1+\rho}}-\frac{1}{\sqrt{1-\rho}}\\[2ex]
\frac{1}{\sqrt{1+\rho}}-\frac{1}{\sqrt{1-\rho}}&\frac{1}{\sqrt{1+\rho}}+\frac{1}{\sqrt{1-\rho}}\\
\end{array}\right].
\end{gather*}
Hence the density of the martingale measure \eqref{gestosc miary martyngalowej} can be written as
\begin{gather}\label{wzor koncowy na tildeZ_T}
\tilde{Z}_T=e^{-A_1W^1_T-A_2W^2_T-BT}=e^{-A_1\widetilde{W}^1_T-A_2\widetilde{W}^2_T-\tilde{B}T},
\end{gather}
where
\begin{align*}
A_1&:=\frac{1}{\rho^2-1}\left(-\frac{\alpha_1-r}{\sigma_1}+\rho\frac{\alpha_2-r}{\sigma_2}\right)\\[1ex]
A_2&:=\frac{1}{\rho^2-1}\left(\rho\frac{\alpha_1-r}{\sigma_1}-\frac{\alpha_2-r}{\sigma_2}\right)\\[1ex]
B&:=\frac{1}{8}\Bigg(\bigg(\Big(\frac{1}{\sqrt{1+\rho}}+\frac{1}{\sqrt{1-\rho}}\Big)\frac{\alpha_1-r}{\sigma_1}
+\Big(\frac{1}{\sqrt{1+\rho}}-\frac{1}{\sqrt{1-\rho}}\Big)\frac{\alpha_2-r}{\sigma_2}\bigg)^2\\[1ex]
&+\bigg(\Big(\frac{1}{\sqrt{1+\rho}}-\frac{1}{\sqrt{1-\rho}}\Big)\frac{\alpha_1-r}{\sigma_1}
+\Big(\frac{1}{\sqrt{1+\rho}}+\frac{1}{\sqrt{1-\rho}}\Big)\frac{\alpha_2-r}{\sigma_2}\bigg)^2
\Bigg)\\[1ex]
\tilde{B}&:=B-A_1\frac{\alpha_1-r}{\sigma_1}-A_2\frac{\alpha_2-r}{\sigma_2}.
\end{align*}

In the following subsections we will use the universal constants:
$A_1, A_2, B, \tilde{B}$ defined in \eqref{wzor koncowy na tildeZ_T} as well as
$a_1, a_2, b, \tilde{a}_1, \tilde{a}_2,
\tilde{b}$ introduced below.\\
\noindent Fix numbers $K>0$, $c\geq 0$. One can check the following
\begin{align}\label{wzor S^1>K}
\left\{S^1_T\geq K\right\}&=\left\{W^1_T\geq a_1\right\}=\left\{\widetilde{W}^1_T\geq \tilde{a}_1\right\},\\[2ex]\label{wzor S^2>K}
\left\{S^2_T\geq K\right\}&=\left\{W^2_T\geq a_2\right\}=\left\{\widetilde{W}^2_T\geq \tilde{a}_2\right\},\\[2ex]\label{wzor S^1>S^2}
\left\{S^1_T\geq
S^2_T\right\}&=\left\{\sigma_1W^1_T-\sigma_2W^2_T\geq
b\right\}=\left\{\sigma_1\widetilde{W}^1_T-\sigma_2\widetilde{W}^2_T\geq
\tilde{b}\right\},\\[2ex]\label{wzor Z^-1>c}
\left\{\tilde{Z}^{-1}_{T}\geq c\right\}&=\left\{A_1 W^{1}_{T}+A_2 W^{2}_{T}\geq \ln c-BT\right\}\\
&=\left\{A_1 \widetilde{W}^{1}_{T}+A_2 \widetilde{W}^{2}_{T}\geq \ln c-\tilde{B}T\right\},
\end{align}
where
\begin{align*}
a_1&:=\frac{1}{\sigma_1}\left(\ln\frac{K}{S^1_0}-(\alpha_1-\frac{1}{2}\sigma_1^2)T\right),\quad
\tilde{a}_1:=\frac{1}{\sigma_1}\left(\ln\frac{K}{S^1_0}-(r-\frac{1}{2}\sigma_1^2)T\right),
\\[2ex]
a_2&:=\frac{1}{\sigma_2}\left(\ln\frac{K}{S^2_0}-(\alpha_2-\frac{1}{2}\sigma_2^2)T\right),\quad
\tilde{a}_2:=\frac{1}{\sigma_2}\left(\ln\frac{K}{S^2_0}-(r-\frac{1}{2}\sigma_2^2)T\right),\\[2ex]
b&:=\ln\left(\frac{S^2_0}{S^1_0}\right)
+(\alpha_2-\alpha_1-\frac{1}{2}(\sigma_2^2-\sigma_1^2))T, \quad
\tilde{b}:=\ln\left(\frac{S^2_0}{S^1_0}\right)
-\frac{1}{2}(\sigma_2^2-\sigma_1^2)T.
\end{align*}
In all the formulas appearing in the sequel it is understood that
$\ln0=-\infty$ and $\Phi$ stands for the distribution function of
$N(0,1)$.

For each derivative below we calculate the risk functions for the case when
when $l(x)=x$ and $l(x)=\frac{x^p}{p},p>1$. In the latter case we use the notation $\Psi^p_1=\Psi^l_1,
\Psi^p_2=\Psi^l_2$. For $l(x)=\frac{x^p}{p}$ we have $I(x)=x^{\frac{1}{p-1}}$ and in
view of \eqref{postac varphi_c} the following holds
\begin{align}\label{wzor ogolny psi-1-l}
\Psi^p_1(c)&=\frac{1}{p}\mathbf{E}\Big[H^{p}\mathbf{1}_{A_c^c}\Big]+\frac{1}{p}\mathbf{E}\Big[(c\tilde{Z}_T)^{\frac{p}{p-1}}\mathbf{1}_{A_c}\Big],\\[2ex]\label{wzor ogolny psi-2-l}
\Psi^p_2(c)&=\mathbf{\tilde{E}}\Big[\Big(H-\big(c\tilde{Z}_T\big)^{\frac{1}{p-1}}\Big)\mathbf{1}_{A_c}\Big].
\end{align}
where
\begin{gather}\label{postac A-c}
 A_c:=\{c\tilde{Z}_T\leq H^{p-1}\},
\end{gather}
and $A_c^c$ stands for the compliment of $A_c$.

\subsection{Digital option}
Digital option is a contract with the payoff function of the form
\begin{gather}\label{digital option}
H=K \cdot \mathbf{1}_{\{S^1_T\geq S^2_T\}},\quad \text{where} \quad
K>0.
\end{gather}
Let $(X,Y)$, $(\tilde{X},\tilde{Y})$ be random vectors defined by $X:=\sigma_1W^{1}_{T}-\sigma_2W^{2}_{T},Y:=A_1W_{T}^{1}+A_2W_{T}^{2}$,
$\tilde{X}:=\sigma_1 \widetilde{W}^1_T-\sigma_2\widetilde{W}^2_T$, $\tilde{Y}:=A_1\widetilde{W}^1_T+A_2\widetilde{W}^2_T$. They are normally distributed under $P$, resp. $\tilde{P}$ and their parameters are given by \eqref{rozklad AX}.

\vskip2mm
\noindent
{\bf Linear loss function}\\
Using \eqref{wzor S^1>S^2} and \eqref{wzor Z^-1>c} we
obtain
\begin{align*}
 \Psi_1(c)&=K\mathbf{E}(\mathbf{1}_{\{S^1_T\geq S^2_T\}}\mathbf{1}_{\{\tilde{Z}_{T}^{-1}\geq c\}})\\
 &=KP(\sigma_1W^{1}_{T}-\sigma_2W^{2}_{T}\geq b,A_1W_{T}^{1}+A_2W_{T}^{2}\geq \ln c-BT),
\end{align*}
and thus
\begin{gather*}
 \Psi_1(c)=K\int_{b}^{+\infty}\int_{\ln c-BT}^{+\infty}f_{X,Y}(x,y)dydx.
\end{gather*}

\noindent
Analogous computation yields
\begin{align*}
\Psi_2(c)&=K\tilde{P}(\sigma_1\widetilde{W}^{1}_{T}-\sigma_2\widetilde{W}^{2}_{T}\geq \tilde{b},A_1\widetilde{W}_{T}^{1}+A_2\widetilde{W}_{T}^{2}\geq \ln c-\tilde{B}T)\\[1ex]
&=K\int_{\tilde{b}}^{+\infty}\int_{\ln c-\tilde{B}T}^{+\infty}\tilde{f}_{\tilde{X},\tilde{Y}}(x,y)dydx.
\end{align*}

\noindent
{\bf Power loss function}\\ In view of \eqref{wzor
S^1>S^2} and \eqref{wzor koncowy na tildeZ_T} we have
\begin{align}\nonumber
A_c&:=\{c\tilde{Z}_T\leq H^{p-1}\}=\{c\tilde{Z}_T\leq
K^{p-1}\mathbf{1}_{\{\sigma_1 W^1_T-\sigma_2 W^2_T\geq b\}}\}\\ \nonumber
&=\{\sigma_1 W^1_T-\sigma_2 W^2_T\geq b,c\tilde{Z}_T\leq K^{p-1}\}\\[1ex]\label{pierwsza postac zbioru}
&=\{\sigma_1 W^1_T-\sigma_2 W^2_T\geq b, A_1W^1_T+A_2W^2_T\geq \ln\left(\frac{K^{p-1}}{c}\right)-BT\}\\[1ex]\label{druga postac zbioru}
&=\{\sigma_1 \widetilde{W}^1_T-\sigma_2 \widetilde{W}^2_T\geq
\tilde{b}, A_1\widetilde{W}^1_T+A_2\widetilde{W}^2_T\geq
\ln\left(\frac{K^{p-1}}{c}\right)-\tilde{B}T\},
\end{align}
and thus
\begin{align*}
\Psi^p_1(c)=\frac{1}{p}\mathbf{E}[K^{p}\mathbf{1}_{\{\sigma_1
W^1_T-\sigma_2 W^2_T\geq
b\}}\mathbf{1}_{A^c_c}]+\frac{1}{p}c^{\frac{p}{p-1}}\mathbf{E}[\tilde{Z}_T^{\frac{p}{p-1}}\mathbf{1}_{A_c}],
\end{align*}
\begin{align*}
\Psi_2^p(c)=\mathbf{\tilde{E}}[K\mathbf{1}_{\{\sigma_1
\widetilde{W}^1_T-\sigma_2 \widetilde{W}^2_T\geq
\tilde{b}\}}\mathbf{1}_{A_c}]-c^\frac{1}{p-1}\mathbf{\tilde{E}}[\tilde{Z}_T^{\frac{1}{p-1}}\mathbf{1}_{A_c}].
\end{align*}

\noindent
In view of \eqref{pierwsza postac zbioru} and \eqref{druga postac zbioru} we have
\begin{align*}
 \Psi_1^p(c)&=\frac{K^p}{p}P\left(\sigma_1 W^1_T-\sigma_2 W^2_T\geq b,A_1W^1_T+A_2W^2_T< \ln\left(\frac{K^{p-1}}{c}\right)-BT\right)\\
 &\phantom{=}+\frac{1}{p}c^{\frac{p}{p-1}}\mathbf{E}[\tilde{Z}_T^{\frac{p}{p-1}}\mathbf{1}_{A_c}]\\[1ex]
&=\frac{K^p}{p}\int_{b}^{+\infty}\int_{-\infty}^{\ln\left(\frac{K^{p-1}}{c}\right)-BT}f_{X,Y}(x,y)dydx\\[1ex]
&\phantom{=}+\frac{1}{p}c^{\frac{p}{p-1}}
\int_{b}^{+\infty}\int_{\ln\left(\frac{K^{p-1}}{c}\right)-BT}^{+\infty}e^{-\frac{p(y+BT)}{p-1}}f_{X,Y}(x,y)dydx,
\end{align*}
and
\begin{align*}
\Psi_2^p(c)&=K\tilde{P}(A_c)-c^{\frac{1}{p-1}}\mathbf{\tilde{E}}[e^{-A_1\widetilde{W}^1_T-A_2\widetilde{W}^2_T-\tilde{B}T}\mathbf{1}_{A_c}]\\[1ex]
&=K\int_{\tilde{b}}^{+\infty}\int_{\ln\left(\frac{K^{p-1}}{c}\right)-\tilde{B}T}^{+\infty}\tilde{f}_{\tilde{X},\tilde{Y}}(x,y)dydx\\[1ex]
&\phantom{=}-c^{\frac{1}{p-1}}\int_{\tilde{b}}^{+\infty}\int_{\ln\left(\frac{K^{p-1}}{c}\right)-\tilde{B}T}^{+\infty}e^{-\frac{y+\tilde{B}T}{p-1}}\tilde{f}_{\tilde{X},\tilde{Y}}(x,y)dydx.
\end{align*}

\subsection{Quantos}
\subsubsection{Quanto domestic}
The contingent claim is of the form
\begin{gather}\label{quanto domestic}
H=S^2_T(S^1_T-K)^{+}, \quad K>0.
\end{gather}

\noindent
{\bf Linear loss function} \\
Using \eqref{wzor S^1>K} we obtain
\begin{align*}
 &\Psi_1(c)=\mathbf{E}[S^2_T(S^{1}_T-K)^{+}\mathbf{1}_{\{\tilde{Z}^{-1}_T\geq c\}}]=\mathbf{E}\Big[S^2_T(S^{1}_T-K)\mathbf{1}_{\{\tilde{Z}^{-1}_T\geq c\}}\mid S^{1}_{T}>K\Big]\cdot\\
 &\cdot P(S^{1}_{T}>K)=\mathbf{E}\Big[S^2_T(S^{1}_T-K)\mathbf{1}_{\{A_1 W^{1}_{T}+A_2 W^{2}_{T}\geq \ln c-BT\}}\mid W^{1}_{T}>a_1\Big]\\[1ex]
&\cdot P(W^{1}_{T}>a_1)=\int_{a_1}^{+\infty}\mathbf{E}\Big[S^2_{0}e^{(\alpha_2-\frac{1}{2}\sigma_2^2)T+\sigma_2W^{2}_{T}}(S^1_{0}e^{(\alpha_1-\frac{1}{2}\sigma_1^2)T+\sigma_1W^{1}_{T}}-K)\cdot\\
&\qquad\qquad\qquad\qquad\qquad\qquad\qquad \cdot\mathbf{1}_{\{W^2_T\geq\frac{\ln c-BT-A_1W^{1}_{T}}{A_2}\}}\mid W_{T}^{1}=x\Big]f_{W_{T}^{1}}(x)dx\\[1ex]
&=S^2_{0}e^{(\alpha_2-\frac{1}{2}\sigma_2^2)T}\int_{a_1}^{+\infty}\int_{\frac{\ln c-BT-A_1x}{A_2}}^{+\infty}(S^{1}_{0}e^{(\alpha_1-\frac{1}{2}\sigma_1^2)T+\sigma_1x}-K)e^{\sigma_2y}\cdot\\
&\qquad\qquad\qquad\qquad\qquad\qquad\qquad\cdot f_{W_{T}^{2}\mid W^{1}_{T}=x}(y)f_{W_{T}^{1}}(x)dy dx,
\end{align*}
and
\begin{align*}
&\Psi_2(c)=\mathbf{\tilde{E}}[S^2_T(S^{1}_T-K)^{+}\mathbf{1}_{\{\tilde{Z}^{-1}_T\geq c\}}]=\mathbf{\tilde{E}}\Big[S^2_T(S^{1}_T-K)\mathbf{1}_{\{\tilde{Z}^{-1}_T\geq c\}}\mid S^{1}_{T}>K\Big]\cdot\\
&\cdot \tilde{P}(S^{1}_{T}>K)=\mathbf{\tilde{E}}\Big[S^2_T(S^{1}_T-K)\mathbf{1}_{\{A_1 \widetilde{W}^{1}_{T}+A_2 \widetilde{W}^{2}_{T}\geq \ln c-\tilde{B}T\}}\mid \widetilde{W}^{1}_{T}>\tilde{a}_1\Big]\cdot\\
&\cdot\tilde{P}(\widetilde{W}^{1}_{T}>\tilde{a}_1)=\int_{\tilde{a}_1}^{+\infty}\mathbf{\tilde{E}}\Big[S^2_{0}e^{(r-\frac{1}{2}\sigma_2^2)T+\sigma_2\widetilde{W}^{2}_{T}}(S^1_{0}e^{(r-\frac{1}{2}\sigma_1^2)T+\sigma_1\widetilde{W}^{1}_{T}}-K)\cdot\\
&\qquad\qquad\qquad\qquad\qquad\qquad\qquad\cdot\mathbf{1}_{\{\widetilde{W}^2_T\geq\frac{\ln c-\tilde{B}T-A_1\widetilde{W}^{1}_{T}}{A_2}\}}\mid \widetilde{W}_{T}^{1}=x\Big]\tilde{f}_{\widetilde{W}_{T}^{1}}(x)dx\\[2ex]
&=S^2_{0}e^{(r-\frac{1}{2}\sigma_2^2)T}\int_{\tilde{a}_1}^{+\infty}\int_{\frac{\ln c-\tilde{B}T-A_1x}{A_2}}^{+\infty}(S^{1}_{0}e^{(r-\frac{1}{2}\sigma_1^2)T+\sigma_1x}-K)e^{\sigma_2y}\cdot\\
&\qquad\qquad\qquad\qquad\qquad\qquad\qquad\cdot\tilde{f}_{\widetilde{W}_{T}^{2}\mid \widetilde{W}^{1}_{T}=x}(y)dy\tilde{f}_{\widetilde{W}_{T}^{1}}(x)dx.
\end{align*}
{\bf Power loss function}\\
\noindent The set \eqref{postac A-c} is of the form
\begin{align*}
A_c&:=\{c\tilde{Z}_T\leq H^{p-1}\}=\left\{\frac{\left(ce^{-A_1W^1_T-A_2W^2_T-BT}\right)^{\frac{1}{p-1}}}{S^2_0e^{(\alpha_2-\frac{1}{2}\sigma_2^2)T+\sigma_2W^2_T}}\leq(S^1_T-K)^{+}\right\}\\[1ex]
&=\left\{\frac{c^{\frac{1}{p-1}}}{S_0^2}e^{-\frac{A_1}{p-1}W^1_T-(\frac{A_2}{p-1}+\sigma_2)W^2_T-(B+\alpha_2-\frac{1}{2}\sigma_2^2)T}\leq
S^1_T-K, S^1_T\geq K)\right\}.
\end{align*}
For simplicity we assume that $\frac{A_2}{p-1}+\sigma_2>0$. In the
opposite case one has to modify the form of the set $A_c$ and thus
also the integration limits in the formulas below. We obtain
\begin{align*}
A_c&=\left\{W_T^2\geq w(W^1_T),W^1_T\geq
a_1\right\}=\left\{\widetilde{W}_T^2\geq
\tilde{w}(\widetilde{W}^1_T),\widetilde{W}^1_T\geq
\tilde{a}_1\right\},
\end{align*}
where
\begin{align*}
 w(x)&:=\frac{\frac{A_1}{p-1}x+\ln\left(\frac{S^2_0(S^1_0e^{(\alpha_1-\sigma_1^2)T+\sigma_1x}-K)}{c^{\frac{1}{p-1}}}\right)+(B+\alpha_2-\frac{1}{2}\sigma_2^2)T}{-(\frac{A_2}{p-1}+\sigma_2)},\\[1ex]
\tilde{w}(x)&:=\frac{\frac{A_1}{p-1}x+\ln\left(\frac{S^2_0(S^1_0e^{(r-\sigma_1^2)T+\sigma_1x}-K)}{c^{\frac{1}{p-1}}}\right)+(\tilde{B}+\alpha_2-\frac{1}{2}\sigma_2^2)T}{-(\frac{A_2}{p-1}+\sigma_2)}.
\end{align*}
In view of this above, \eqref{wzor ogolny psi-1-l}, \eqref{wzor
ogolny psi-2-l} and using conditional densities we obtain
\begin{align*}
 &\Psi_1^p(c)=\frac{(S_0^2)^pe^{(\alpha_2-\frac{1}{2}\sigma_2^2)pT}}{p}\bigg(\int_{a_1}^{+\infty}\int_{-\infty}^{+\infty}e^{p\sigma_2y}(S_0^1e^{(\alpha_1-\frac{1}{2}
\sigma_1^2)T+\sigma_2x}-K)^p\cdot\\ 
&\qquad\qquad\qquad\qquad\qquad\qquad\qquad\cdot f_{W^2_T\mid W^1_T=x}(y)f_{W^1_T}(x)dydx\\[1ex]
&-\int_{a_1}^{+\infty}\int_{w(x)}^{+\infty}e^{p\sigma_2y}(S_0^1e^{(\alpha_1-\frac{1}{2}
\sigma_1^2)T+\sigma_2x}-K)^pf_{W^2_T\mid W^1_T=x}(y)f_{W^1_T}(x)dydx\bigg)\\[1ex]
&+\frac{c^{\frac{p}{p-1}}e^{-\frac{BTp}{p-1}}}{p}\int_{a_1}^{+\infty}\int_{w(x)}^{+\infty}e^{-\left(\frac{A_1p}{p-1}x+\frac{A_2p}{p-1}y\right)}f_{W^2_T\mid W^1_T=x}(y)f_{W^1_T}(x)dydx,\\
\text{and}\\
 &\Psi_2^p(c)=S_0^2e^{(r-\frac{1}{2}\sigma_2^2)T}\int_{\tilde{a}_1}^{+\infty}\int_{\tilde{w}(x)}^{+\infty}e^{\sigma_2y}(S_0^1e^{(r-\frac{1}{2}
\sigma_1^2)T+\sigma_2x}-K)\cdot\\
&\qquad\qquad\qquad\qquad\qquad\qquad\qquad\cdot\tilde{f}_{\widetilde{W}^2_T\mid\widetilde{W}^1_T=x}(y)\tilde{f}_{\widetilde{W}^1_T}(x)dydx\\[1ex]
&-c^{\frac{1}{p-1}}e^{-\frac{BT}{p-1}}\int_{\tilde{a}_1}^{+\infty}\int_{\tilde{w}(x)}^{+\infty}e^{-\frac{1}{p-1}(A_1x+A_2y)}\tilde{f}_{\widetilde{W}^2_T\mid\widetilde{W}^1_T=x}(y)\tilde{f}_{\widetilde{W}^1_T}(x)dydx.
\end{align*}

\subsubsection{Quanto foreign}
The payoff is of the form
\begin{gather*}
H=\left(S^1_T-\frac{K}{S^2_T}\right)^+, \quad K>0.
\end{gather*}
\noindent {\bf Linear loss function}\\
First let us notice that
\begin{gather}\label{wzor 1 QF}
\left\{S^1_T-\frac{K}{S^2_T}\geq0\right\}=\left\{\sigma_1W^1_T+\sigma_2W^2_T\geq
d\right\}=\left\{\sigma_1\widetilde{W}^1_T+\sigma_2\widetilde{W}^2_T\geq
\tilde{d}\right\},
\end{gather}
where
{\small
\begin{align}\label{QF stala d}
d:=\ln\frac{K}{S^1_0S^2_0}-\left(\alpha_1+\alpha_2-\frac{1}{2}(\sigma_1^2+\sigma_2^2)\right)T,\quad
\tilde{d}:=\ln\frac{K}{S^1_0S^2_0}-\left(2r-\frac{1}{2}(\sigma_1^2+\sigma_2^2)\right)T.
\end{align}}
\noindent
We have
\begin{align*}
\Psi_1(c)&=\mathbf{E}\Big[\left(S^1_T-\frac{K}{S^2_T}\right)^+\mathbf{1}_{\{\tilde{Z}^{-1}_T\geq c\}}\Big]\\
&=\mathbf{E}\Big[\left(S^1_T-\frac{K}{S^2_T}\right)\mathbf{1}_{\{W^2_T\geq\frac{\ln c-BT-A_1W^{1}_{T}}{A_2}\}}\mid \sigma_1W^1_T+\sigma_2W^2_T\geq
d\Big]\cdot\\
&\qquad\qquad\qquad\qquad\qquad\qquad\qquad\cdot P(\sigma_1W^1_T+\sigma_2W^2_T\geq d).
\end{align*}
Denoting $Z:=\sigma_1W^1_T+\sigma_2W^2_T$ and taking into account
conditional distribution $\mathcal{L}(W^{1}_T,W^2_T~\mid~Z)$ we
obtain {\small
\begin{align*}
\Psi_1(c)&=\int_{d}^{+\infty}\int_{-\infty}^{+\infty}\int_{\frac{\ln c-BT-A_1x}{A_2}}^{+\infty}(S^1_0e^{(\alpha_1-\frac{1}{2}\sigma_1^2)T+\sigma_1x}-KS_{0}^{2}e^{(-\alpha_2+\frac{1}{2}\sigma_2^2)T-
\sigma_2y})\cdot\\ 
&\qquad\qquad\qquad\qquad\qquad\qquad\qquad\cdot f_{(W^{1}_{T},W^{2}_{T})\mid Z=z}(x,y)dydxf_Z(z)dz.
\end{align*}}
Using the same argument under the measure $\tilde{P}$ with
$\tilde{Z}:=\sigma_1\widetilde{W}^1_T+\sigma_2\widetilde{W}^2_T$
yields {\small
\begin{align*}
 \Psi_2(c)&=\mathbf{\tilde{E}}\Big[\left(S^1_T-\frac{K}{S^2_T}\right)\mathbf{1}_{\{\widetilde{W}^2_T\geq\frac{\ln c-\tilde{B}T-A_1\widetilde{W}^{1}_{T}}{A_2}\}}\mid \sigma_1\widetilde{W}^1_T+\sigma_2\widetilde{W}^2_T\geq
\tilde{d}\Big]\cdot\\
&\qquad\qquad\qquad\qquad\qquad\qquad\qquad\cdot\tilde{P}(\sigma_1\widetilde{W}^1_T+\sigma_2\widetilde{W}^2_T\geq
\tilde{d})\\[2ex]
&=\int_{\tilde{d}}^{+\infty}\int_{-\infty}^{+\infty}\int_{\frac{\ln c-\tilde{B}T-A_1x}{A_2}}^{+\infty}(S^1_0e^{(r-\frac{1}{2}\sigma_1^2)T+\sigma_1x}-KS_{0}^{2}e^{(-r+\frac{1}{2}\sigma_2^2)T-
\sigma_2y})\cdot\\
&\qquad\qquad\qquad\qquad\qquad\qquad\qquad\cdot\tilde{f}_{(\widetilde{W}^{1}_{T},\widetilde{W}^{2}_{T})\mid  \tilde{Z}=z}(x,y)dydx\tilde{f}_{\tilde{Z}}(z)dz.
\end{align*}}
\noindent
{\bf Power loss function}\\
Using \eqref{wzor 1 QF} one can check the following
\begin{align}\nonumber
A_c&:=\Big\{c\tilde{Z}_T\leq
\Big(\Big(S^1_T-\frac{K}{S^2_T}\Big)^{+}\Big)^{p-1},S^1_T-\frac{K}{S^2_T}>0\Big\}\\[1ex]\nonumber
&=\Big\{c\tilde{Z}_T\leq
\Big(\Big(S^1_T-\frac{K}{S^2_T}\Big)^{+}\Big)^{p-1},\sigma_1W^{1}_{T}+\sigma_2W^{2}_{T}>d\Big\}\\[1ex]\label{QF zbior Ac
1}
&=\Big\{\frac{A_1}{p-1}W^{1}_{T}+\big(\frac{A_2}{p-1}-\sigma_2\big)W^2_T\geq
v(\sigma_1W^{1}_{T}+\sigma_2W^{2}_{T}),\sigma_1W^{1}_{T}+\sigma_2W^{2}_{T}>d\Big\}\\[1ex]\label{QF zbior Ac 2}
&=\Big\{\frac{A_1}{p-1}\widetilde{W}^{1}_{T}+\big(\frac{A_2}{p-1}-\sigma_2\big)\widetilde{W}^2_T\geq
\tilde{v}(\sigma_1\widetilde{W}^{1}_{T}+\sigma_2\widetilde{W}^{2}_{T}),\sigma_1\widetilde{W}^{1}_{T}+\sigma_2\widetilde{W}^{2}_{T}>\tilde{d}\Big\},
\end{align}
where $d,\tilde{d}$ are given by \eqref{QF stala d} and
\begin{align*}
v(x)&=\ln\left\{\frac{S^1_0S^2_0e^{(\alpha_1+\alpha_2-\frac{1}{2}(\sigma_1^2+\sigma_2^2))T+x}-K}{c^{\frac{1}{p-1}}S^2_0e^{(\alpha_2-\frac{1}{2}\sigma_2^2-\frac{B}{p-1})T}}\right\},\\[2ex]
\tilde{v}(x)&=\ln\left\{\frac{S^1_0S^2_0e^{(2r-\frac{1}{2}(\sigma_1^2+\sigma_2^2))T+x}-K}{c^{\frac{1}{p-1}}S^2_0e^{(r-\frac{1}{2}\sigma_2^2-\frac{\tilde{B}}{p-1})T}}\right\}.
\end{align*}

\noindent To calculate $\Psi_1^l,\Psi_2^l$ we use conditional
distributions $\mathcal{L}(X\mid Y)$,
$\mathcal{L}(\tilde{X}\mid\tilde{Y})$, where
$X:=\frac{A_1}{p-1}W^{1}_{T}+\big(\frac{A_2}{p-1}-\sigma_2\big)W^2_T$,
$Y:=\sigma_1W^1_T+\sigma_2W^2_T$,
$\tilde{Y}:=\frac{A_1}{p-1}\widetilde{W}^{1}_{T}+\big(\frac{A_2}{p-1}-\sigma_2\big)\widetilde{W}^2_T$,
$\tilde{Y}:=\sigma_1\widetilde{W}^1_T+\sigma_2\widetilde{W}^2_T$.
Denote by $k_1,k_2,k_3,k_4$ constants satisfying $W^1_T=k_1X+k_2Y$,
$W^2_T=k_3X+k_4Y$,
$\widetilde{W}^1_T=k_1\widetilde{X}+k_2\widetilde{Y}$,
$\widetilde{W}^2_T=k_3\widetilde{X}+k_4\widetilde{Y}$. Then we have
{\small
\begin{align*}
&\Psi_1^p(c)=\frac{1}{p}\int_{d}^{\infty}\int^{v(y)}_{-\infty}\left(S^1_0e^{(\alpha_1-\frac{1}{2}\sigma_1^2)T+\sigma_1(k_1x+k_2y)}-\frac{K}{S^1_0e^{(\alpha_2-\frac{1}{2}\sigma_2^2)T+\sigma_1(k_3x+k_4y)}}\right)^p
\cdot\\[1ex]
&\qquad\qquad\qquad\qquad\qquad\qquad\qquad\qquad\qquad\cdot{f}_{X\mid Y=y}(x)f_{Y}(y)dxdy\\[1ex]
&+\frac{1}{p}c^{\frac{p}{p-1}}e^{-\frac{pBT}{p-1}}\int_{d}^{+\infty}\int_{v(y)}^{+\infty}
e^{-\frac{pA_1}{p-1}(k_1x+k_2y)-\frac{pA_2}{p-1}(k_3x+k_4y)}
f_{X\mid Y=y}(x)f_{Y}(y)dxdy,\\
\text{and}\\
&\Psi_2^p(c)=\int_{\tilde{d}}^{+\infty}\int_{\tilde{v}(y)}^{+\infty}\left(S^1_0e^{(r-\frac{1}{2}\sigma_1^2)T+\sigma_1(k_1x+k_2y)}-\frac{K}{S^1_0e^{(r-\frac{1}{2}\sigma_2^2)T+\sigma_1(k_3x+k_4y)}}\right)\cdot\\[1ex]
&\qquad\qquad\qquad\qquad\qquad\qquad\qquad\qquad\qquad\cdot\tilde{f}_{\tilde{X}\mid
\tilde{Y}=y}(x)\tilde{f}_{\tilde{Y}}(y)dxdy\\[1ex]
&-c^{\frac{1}{p-1}}e^{-\frac{\tilde{B}T}{p-1}}\int_{\tilde{d}}^{+\infty}\int_{\tilde{v}(y)}^{+\infty}
e^{-\frac{A_1}{p-1}(k_1x+k_2y)-\frac{A_2}{p-1}(k_3x+k_4y)}
\tilde{f}_{\tilde{X}\mid
\tilde{Y}=y}(x)\tilde{f}_{\tilde{Y}}(y)dxdy.
\end{align*}}

\subsection{Outperformance option}
The problem is studied for
\begin{gather*}
H=\left(\max\{S^1_T,S^2_T\}-K\right)^{+}, \quad K>0.
\end{gather*}

\noindent
{\bf Linear loss function}\\
By \eqref{wzor S^1>K},
\eqref{wzor S^2>K} and \eqref{wzor S^1>S^2} we get
\begin{align*}
 &\Psi_1(c)=\mathbf{E}\Big[(S_{T}^{1}-K)\mathbf{1}_{\{\tilde{Z}^{-1}_T\geq c\}}\mid S^1_T\geq K, S^1_T\geq S^2_T\Big]P(S^1_T\geq K, S^1_T\geq S^2_T)\\[1ex]
&+\mathbf{E}\Big[(S_{T}^{2}-K)\mathbf{1}_{\{\tilde{Z}^{-1}_T\geq c\}}\mid S^2_T\geq K, S^1_T<S^2_T\Big]P(S^2_T\geq K, S^1_T<S^2_T)\\[1ex]
&=\mathbf{E}\Big[(S_{T}^{1}-K)\mathbf{1}_{\{\tilde{Z}^{-1}_T\geq c\}}\mid W^1_T\geq a_1, \sigma_1W^1_T-\sigma_2W^2_T\geq b\Big]\cdot\\
&\qquad\qquad\qquad\qquad\qquad\qquad\qquad\cdot P(W^1_T\geq a_1, \sigma_1W^1_T-\sigma_2W^2_T\geq b)\\[1ex]
&+\mathbf{E}\Big[(S_{T}^{2}-K)\mathbf{1}_{\{\tilde{Z}^{-1}_T\geq c\}}\mid W^2_T\geq a_2, \sigma_1W^1_T-\sigma_2W^2_T< b\Big]\cdot\\
&\qquad\qquad\qquad\qquad\qquad\qquad\qquad\cdot P(W^2_T\geq a_2, \sigma_1W^1_T-\sigma_2W^2_T< b),
\end{align*}
and further
\begin{align*}
&\Psi_1(c)=\int_{a_1}^{+\infty}\int_{b}^{+\infty}(S^1_0e^{(\alpha_1-\frac{1}{2}\sigma_1^2)T+\sigma_1x}-K)\mathbf{1}_{\{A_1x+A_2\frac{\sigma_1x-z}{\sigma_2}\geq\ln c -BT\}}\cdot\\
&\qquad\qquad\qquad\qquad\qquad\qquad\qquad\qquad\qquad\cdot f_{W^1_T,\sigma_1W^1_T-\sigma_2W^2_T}(x,z)dzdx\\[1ex]
&+\int_{a_2}^{+\infty}\int_{-\infty}^{b}(S^2_0e^{(\alpha_2-\frac{1}{2}\sigma_2^2)T+\sigma_2y}-K)\mathbf{1}_{\{A_1\frac{z+\sigma_2y}{\sigma_1}+A_2y\geq\ln c -BT\}}\cdot\\ 
&\qquad\qquad\qquad\qquad\qquad\qquad\qquad\qquad\qquad\cdot f_{W^2_T,\sigma_1W^1_T-\sigma_2W^2_T}(y,z)dzdy.
\end{align*}
Similarly, for $\Psi_2$ we get 
\begin{align*}
 \Psi_2(c)&=\mathbf{\tilde{E}}\Big[(S_{T}^{1}-K)\mathbf{1}_{\{\tilde{Z}^{-1}_T\geq c\}}\mid S^1_T\geq K, S^1_T\geq S^2_T\Big]\tilde{P}(S^1_T\geq K, S^1_T\geq S^2_T)\\[1ex]
&+\mathbf{\tilde{E}}\Big[(S_{T}^{2}-K)\mathbf{1}_{\{\tilde{Z}^{-1}_T\geq c\}}\mid S^2_T\geq K, S^1_T<S^2_T\Big]\tilde{P}(S^2_T\geq K, S^1_T<S^2_T)\\[1ex]
&=\mathbf{\tilde{E}}\Big[(S_{T}^{1}-K)\mathbf{1}_{\{\tilde{Z}^{-1}_T\geq c\}}\mid \widetilde{W}^1_T\geq \tilde{a}_1, \sigma_1\widetilde{W}^1_T-\sigma_2\widetilde{W}^2_T\geq \tilde{b}\Big]\cdot\\
&\qquad\qquad\qquad\qquad\qquad\qquad\qquad\cdot\tilde{P}(\widetilde{W}^1_T\geq \tilde{a}_1, \sigma_1\widetilde{W}^1_T-\sigma_2\widetilde{W}^2_T\geq \tilde{b})\\[1ex]
&+\mathbf{\tilde{E}}\Big[(S_{T}^{2}-K)\mathbf{1}_{\{\tilde{Z}^{-1}_T\geq c\}}\mid \widetilde{W}^2_T\geq \tilde{a}_2, \sigma_1\widetilde{W}^1_T-\sigma_2\widetilde{W}^2_T< \tilde{b}\Big]\cdot\\
&\qquad\qquad\qquad\qquad\qquad\qquad\qquad\cdot\tilde{P}(\widetilde{W}^2_T\geq \tilde{a}_2, \sigma_1\widetilde{W}^1_T-\sigma_2\widetilde{W}^2_T< \tilde{b})\\[1ex]
\end{align*}
which leads to
\begin{align*}
\Psi_2(c)&=\int_{\tilde{a}_1}^{+\infty}\int_{\tilde{b}}^{+\infty}(S^1_0e^{(r-\frac{1}{2}\sigma_1^2)T+\sigma_1x}-K)\mathbf{1}_{\{A_1x+A_2\frac{\sigma_1x-z}{\sigma_2}\geq\ln c -\tilde{B}T\}}\cdot\\
&\qquad\qquad\qquad\qquad\qquad\qquad\qquad\qquad\qquad\cdot\tilde{f}_{\widetilde{W}^1_T,\sigma_1\widetilde{W}^1_T-\sigma_2\widetilde{W}^2_T}(x,z)dzdx\\[1ex]
&+\int_{\tilde{a}_2}^{+\infty}\int_{-\infty}^{\tilde{b}}(S^2_0e^{(r-\frac{1}{2}\sigma_2^2)T+\sigma_2y}-K)\mathbf{1}_{\{A_1\frac{z+\sigma_2y}{\sigma_1}+A_2y\geq\ln c -\tilde{B}T\}}\cdot\\
&\qquad\qquad\qquad\qquad\qquad\qquad\qquad\qquad\qquad\cdot\tilde{f}_{\widetilde{W}^2_T,\sigma_1\widetilde{W}^1_T-\sigma_2\widetilde{W}^2_T}(y,z)dzdy.
\end{align*}

\noindent
{\bf Power loss function}\\
Taking into account \eqref{wzor S^1>K}, \eqref{wzor S^2>K},
\eqref{wzor S^1>S^2} we can write
\begin{align*}
 A_c=\{&c\tilde{Z}_T\leq(S^1_T\vee S^2_T-K)^{p-1},S^1_T\vee S^2_T-K>0\}=\{c\tilde{Z}_T\leq(S^1_T-K)^{p-1},\\
 & S^1_T>K, S^1_T\geq S^2_T\}\cup
 \{c\tilde{Z}_T\leq(S^2_T-K)^{p-1},S^2_T>K, S^1_T\leq S^2_T\}.
 \end{align*}
We consider the case when $A_1>0,A_2>0$:
  \begin{align}\label{Outperformace A_c l}\nonumber
 &A_c=\{W^2_T\geq-\Big(A_1W^1_T+BT+\ln\Big(\frac{1}{c}(S^1_T-K)^{p-1})\Big)\Big),W^1_T>a_1,\\ \nonumber
 &\sigma_1W^1_T-\sigma_2W^2_T\geq b\}\cup\{W^1_T\geq-\Big(A_2W^2_T+BT+\ln\Big(\frac{1}{c}(S^2_T-K)^{p-1})\Big)\Big),\\ \nonumber &W^2_T>a_2,\sigma_1W^1_T-\sigma_2W^2_T\leq b\}=\left\{W^2_T\geq v_1(W^1_T),W^1_T>a_1,W^2_T\leq\frac{\sigma_1 W^1_T-b}{\sigma_2}\right\}\\[1ex]\nonumber
 &\cup \left\{W^1_T\geq v_2(W^2_T),W^2_T>a_2,W^1_T\leq\frac{\sigma_2 W^2_T-b}{\sigma_1}\right\}\\[1ex]\nonumber
 &=\left\{\widetilde{W}^2_T\geq \tilde{v}_1(\widetilde{W}^1_T),\widetilde{W}^1_T>\tilde{a}_1,\widetilde{W}^2_T\leq\frac{\sigma_1 \widetilde{W}^1_T-\tilde{b}}{\sigma_2}\right\}\\[1ex]
 &\cup \left\{\widetilde{W}^1_T\geq \tilde{v}_2(\widetilde{W}^2_T),\widetilde{W}^2_T>\tilde{a}_2,\widetilde{W}^1_T\leq\frac{\sigma_2 \widetilde{W}^2_T-\tilde{b}}{\sigma_1}\right\},
 \end{align}
where
\begin{align*}
 v_1(x)&=-\frac{1}{A_2}\left(A_1x+BT+\ln\Big(\frac{1}{c}(S^1_0e^{(\alpha_1-\frac{1}{2}\sigma_1^2)T+\sigma_1 x}-K)^p\Big)\right),\\[1ex]
 v_2(x)&=-\frac{1}{A_1}\left(A_2x+BT+\ln\Big(\frac{1}{c}(S^2_0e^{(\alpha_2-\frac{1}{2}\sigma_2^2)T+\sigma_2 x}-K)^p\Big)\right),\\[1ex]
 \tilde{v}_1(x)&=-\frac{1}{A_2}\left(A_1x+\tilde{B}T+\ln\Big(\frac{1}{c}(S^1_0e^{(r-\frac{1}{2}\sigma_1^2)T+\sigma_1 x}-K)^p\Big)\right),\\[1ex]
 \tilde{v}_2(x)&=-\frac{1}{A_1}\left(A_2x+\tilde{B}T+\ln\Big(\frac{1}{c}(S^2_0e^{(r-\frac{1}{2}\sigma_2^2)T+\sigma_2 x}-K)^p\Big)\right).
\end{align*}
Using the representation \eqref{Outperformace A_c l} and accepting
the convention that the integral over the empty set is zero, we
obtain

\begin{align*}
 \Psi^p_1(c)&=\int_{a_1}^{+\infty}\int_{-\infty}^{v_1(x)\wedge\frac{\sigma_1x-b}{\sigma_2}}
 (S^1_0e^{(\alpha_1-\frac{1}{2}\sigma_1^2)T+\sigma_1 x}-K)^p
 {f}_{{W}^2_T\mid {W}^1_T=x}(y)dy f_{W^1_T}(x)dx\\[1ex]
  &+\frac{1}{p}c^{\frac{p}{p-1}}e^{-\frac{pB}{p-1}}\int_{a_1}^{+\infty}\int_{v_1(x)}^{\frac{\sigma_1x-b}{\sigma_2}}
 e^{-\frac{pA_1}{p-1}x-\frac{pA_2}{p-1}y}
 {f}_{{W}^2_T\mid {W}^1_T=x}(y)dy f_{W^1_T}(x)dx\\[1ex]
 &+\int_{a_2}^{+\infty}\int_{-\infty}^{v_2(x)\wedge\frac{\sigma_2x-b}{\sigma_1}}
 (S^2_0e^{(\alpha_2-\frac{1}{2}\sigma_2^2)T+\sigma_2 x}-K)^p
 {f}_{{W}^1_T\mid {W}^2_T=x}(y)dy f_{W^2_T}(x)dx\\[1ex]
 &+\frac{1}{p}c^{\frac{p}{p-1}}e^{-\frac{pB}{p-1}}\int_{a_2}^{+\infty}\int_{v_2(x)}^{\frac{\sigma_2x-b}{\sigma_1}}
 e^{-\frac{pA_1}{p-1}x-\frac{pA_2}{p-1}y}
 {f}_{{W}^1_T\mid {W}^2_T=x}(y)dy f_{W^2_T}(x)dx,
\end{align*}

\begin{align*}
 \Psi^p_2(c)&=\int_{\tilde{a}_1}^{+\infty}\int_{\tilde{v}_1(x)}^{\frac{\sigma_1x-\tilde{b}}{\sigma_2}}(S^1_0e^{(r-\frac{1}{2}\sigma_1^2)T+\sigma_1 x}-K)\tilde{f}_{\widetilde{W}^2_T\mid\widetilde{W}^1_T=x}(y)dy\tilde{f}_{\widetilde{W}^1_T}(x)dx\\[1ex]
 &-c^{\frac{1}{p-1}}e^{-\frac{\tilde{B}}{p-1}T}\int_{\tilde{a}_1}^{+\infty}\int_{\tilde{v}_1(x)}^{\frac{\sigma_1x-\tilde{b}}{\sigma_2}}
 e^{-\frac{A_1}{p-1}x-\frac{A_2}{p-1}y}
 \tilde{f}_{\widetilde{W}^2_T\mid\widetilde{W}^1_T=x}(y)dy\tilde{f}_{\widetilde{W}^1_T}(x)dx\\[1ex]
 &+\int_{\tilde{a}_2}^{+\infty}\int_{\tilde{v}_2(x)}^{\frac{\sigma_2x-\tilde{b}}{\sigma_1}}(S^2_0e^{(r-\frac{1}{2}\sigma_2^2)T+\sigma_2 x}-K)\tilde{f}_{\widetilde{W}^1_T\mid\widetilde{W}^2_T=x}(y)dy\tilde{f}_{\widetilde{W}^2_T}(x)dx\\[1ex]
 &-c^{\frac{1}{p-1}}e^{-\frac{\tilde{B}}{p-1}T}\int_{\tilde{a}_2}^{+\infty}\int_{\tilde{v}_2(x)}^{\frac{\sigma_1x-\tilde{b}}{\sigma_1}}
 e^{-\frac{A_1}{p-1}x-\frac{A_2}{p-1}y}
 \tilde{f}_{\widetilde{W}^1_T\mid\widetilde{W}^2_T=x}(y)dy\tilde{f}_{\widetilde{W}^2_T}(x)dx.
\end{align*}

\subsection{Spread option}

The payoff is of the form
\begin{gather*}
H=\left(S^1_T-S^2_T-K\right)^{+}, \quad K>0.
\end{gather*}
One can check the following
\begin{gather*}
\{S^1_T\geq S^2_T+K\}=\{W^1_T\geq d(W^2_T)\}=\{\widetilde{W}^1_T\geq \tilde{d}(\widetilde{W}^2_T)\},
\end{gather*}
where
\begin{gather*}
d(y):=\frac{1}{\sigma_1}\ln\left(\frac{S^2_0e^{(\alpha_2-\frac{1}{2}\sigma_2^2)T+\sigma_2y}+K}{S^1_0e^{(\alpha_1-\frac{1}{2}\sigma_1^2)T}}\right),\ \tilde{d}(y):=\frac{1}{\sigma_1}\ln\left(\frac{S^2_0e^{(r-\frac{1}{2}\sigma_2^2)T+\sigma_2y}+K}{S^1_0e^{(r-\frac{1}{2}\sigma_1^2)T}}\right).
\end{gather*}
{\bf Linear loss function} \\
We have {\small\begin{align*}
&\Psi_1(c)=\mathbf{E}\Big[(S^1_T-S^2_T-K)^{+}\mathbf{1}_{\{\tilde{Z}^{-1}_T\geq
c\}}\Big]=
\int_{-\infty}^{+\infty}\mathbf{E}\big[(S^1_T-S^2_T-K)^{+}\mathbf{1}_{\{\tilde{Z}^{-1}_T\geq c\}}\mid W^2_T=y\big]\cdot\\
&\qquad\qquad\cdot f_{W^2_T}(y)dy=\int_{-\infty}^{+\infty}\int_{d(y)}^{+\infty}\Big(S^1_0e^{(\alpha_1-\frac{1}{2}\sigma_1^2)T+\sigma_1x}-S^2_0e^{(\alpha_2-\frac{1}{2}\sigma_2^2)T+\sigma_2y}-K\Big)\cdot\\
&\qquad\qquad\cdot \mathbf{1}_{\{A_1x+A_2y\geq\ln c-BT\}}f_{W^1_T\mid W^2_T=y}(x)dxf_{W^2_T}(y)dy,
\end{align*}}
and {\small
\begin{align*}
&\Psi_2(c)=\mathbf{\tilde{E}}\Big[(S^1_T-S^2_T-K)^{+}\mathbf{1}_{\{\tilde{Z}^{-1}_T\geq c\}}\Big]=
\int_{-\infty}^{+\infty}\mathbf{\tilde{E}}\big[(S^1_T-S^2_T-K)^{+}\mathbf{1}_{\{\tilde{Z}^{-1}_T\geq c\}}\mid \widetilde{W}^2_T=y\big]\cdot\\
&\qquad\qquad\cdot\tilde{f}_{\widetilde{W}^2_T}(y)dy\\[1ex]
&=\int_{-\infty}^{+\infty}\int_{\tilde{d}(y)}^{+\infty}\Big(S^1_0e^{(r-\frac{1}{2}\sigma_1^2)T+\sigma_1x}-S^2_0e^{(r-\frac{1}{2}\sigma_2^2)T+\sigma_2y}-K\Big)\\
&\qquad\qquad\cdot\mathbf{1}_{\{A_1x+A_2y\geq\ln c-\tilde{B}T\}}\tilde{f}_{\widetilde{W}^1_T\mid \widetilde{W}^2_T=y}(x)dxf_{W^2_T}(y)dy.
\end{align*}}

\noindent
{\bf Power loss function} \\
We have {\small
\begin{align}\label{spread option zbior A_c}\nonumber
 &A_c:=\{c\tilde{Z}_T\leq(S^1_T-S^2_T-K)^{p-1},S^1_T-S^2_T-K>0\}\\[1ex]\nonumber
 &=\Big\{c^{\frac{1}{p-1}}e^{-\frac{A_1}{p-1}W^1_T-\frac{A_2}{p-1}W^2_T-\frac{B}{p-1}T}\leq S^1_0e^{(\alpha_1-\frac{1}{2}\sigma_1^2)T+\sigma_1W^1_T}-
 S^2_0e^{(\alpha_2-\frac{1}{2}\sigma^2_2)T+\sigma_2W^2_T}-K,\\ 
&\qquad W^1_T\geq d(W^2_T)\Big\}=\Big\{W^1_T\in \mathcal{A}(W^2_T)\Big\}=\Big\{\widetilde{W}^1_T\in \tilde{\mathcal{A}}(\widetilde{W}^2_T)\Big\},
\end{align}}
where
\begin{align*}
\mathcal{A}(y):=\{&x:c^{\frac{1}{p-1}}e^{-\frac{A_1}{p-1}x-\frac{A_2}{p-1}y-\frac{B}{p-1}T}\leq S^1_0e^{(\alpha_1-\frac{1}{2}\sigma_1^2)T+\sigma_1x}-
 S^2_0e^{(\alpha_2-\frac{1}{2}\sigma^2_2)T+\sigma_2y}-K,\\ 
 &x\geq d(y)\},\\
\tilde{\mathcal{A}}(y):=\{&x:c^{\frac{1}{p-1}}e^{-\frac{A_1}{p-1}x-\frac{A_2}{p-1}y-\frac{\tilde{B}}{p-1}T}\leq S^1_0e^{(r-\frac{1}{2}\sigma_1^2)T+\sigma_1x}-
 S^2_0e^{(r-\frac{1}{2}\sigma^2_2)T+\sigma_2y}-K,\\ 
 & x\geq \tilde{d}(y)\}.
\end{align*}
Let us notice that the set $A^c_c\cap\{H>0\}$ is of the form
\begin{gather}\label{spread option postac A_c^c}
 A^c_c\cap\{H>0\}=\{W^1_T\in\mathcal{B}(W^2_T)\},
\end{gather}
where
\begin{align*}
\mathcal{B}(y):=\{&x:c^{\frac{1}{p-1}}e^{-\frac{A_1}{p-1}x-\frac{A_2}{p-1}y-\frac{B}{p-1}T}> S^1_0e^{(\alpha_1-\frac{1}{2}\sigma_1^2)T+\sigma_1x}-
 S^2_0e^{(\alpha_2-\frac{1}{2}\sigma^2_2)T+\sigma_2y}-K,\\ 
 &x\geq d(y)\}.
\end{align*}
Taking into account \eqref{spread option zbior A_c} and \eqref{spread option postac A_c^c} we obtain
\begin{align*}
 \Psi_1^p(c)&=
\frac{1}{p}\int_{-\infty}^{+\infty}\int_{\mathcal{B}(y)}\Big(
S^1_0e^{(\alpha_1-\frac{1}{2}\sigma_1^2)T+\sigma_1x}-
 S^2_0e^{(\alpha_2-\frac{1}{2}\sigma^2_2)T+\sigma_2y}-K
\Big)^p\cdot\\
&\qquad\qquad\qquad\qquad\qquad\qquad\cdot f_{W^1_T\mid W^2_T=y}(x)dx f_{W^2_T}(y)dy\\[1ex]
 &+\frac{1}{p}c^{\frac{p}{p-1}}e^{-\frac{pBT}{p-1}}\int_{-\infty}^{+\infty}\int_{\mathcal{A}(y)}\Big(e^{-\frac{pA_1}{p-1}x-\frac{pA_2}{p-1}y}\Big)f_{W^1_T\mid W^2_T=y}(x)dx f_{W^2_T}(y)dy,
\end{align*}

\begin{align*}
 \Psi_2^p(c)&=\int_{-\infty}^{+\infty}\int_{\tilde{\mathcal{A}}(y)}\Big(S^1_0e^{(r-\frac{1}{2}\sigma_1^2)T+\sigma_1x}-
 S^2_0e^{(r-\frac{1}{2}\sigma^2_2)T+\sigma_2y}-K\Big)\cdot\\
 &\qquad\qquad\qquad\qquad\qquad\qquad\cdot\tilde{f}_{\widetilde{W}^1_T\mid\widetilde{W}^2_T=y}(x)dx\tilde{f}_{\widetilde{W}^2_T}(y)dy\\[1ex]
 &+c^{\frac{1}{p-1}}e^{-\frac{\tilde{B}T}{p-1}}\int_{-\infty}^{+\infty}\int_{\tilde{\mathcal{A}}(y)}\Big(e^{-\frac{A_1}{p-1}x-\frac{A_2}{p-1}y}\Big)\tilde{f}_{\widetilde{W}^1_T\mid\widetilde{W}^2_T=y}(x)dx\tilde{f}_{\widetilde{W}^2_T}(y)dy.
\end{align*}

\subsection*{Acknowledgements}
Research supported by the Polish MNiSW grant NN201419039.

\end{document}